\documentclass[11pt]{amsart}
\usepackage{amsxtra}
\usepackage{amssymb,mathrsfs}
\usepackage{mathtools}
\usepackage{color}
\usepackage{tikz}
\usetikzlibrary{cd}
\usetikzlibrary{arrows}
\usetikzlibrary{matrix}
\theoremstyle{plain}

\newcommand{\cleqn}{\setcounter{equation}{0}}
\newcommand{\clth}{\setcounter{theorem}{0}}
\newcommand {\sectionnew}[1]{\section{#1}\cleqn\clth}
\newcommand{\nn}{\hfill\nonumber}
\newtheorem{theorem}{Theorem}[section]
\newtheorem{lemma}[theorem]{Lemma}
\newtheorem{definition-theorem}[theorem]{Definition-Theorem}
\newtheorem{proposition}[theorem]{Proposition}
\newtheorem{corollary}[theorem]{Corollary}
\newtheorem{definition}[theorem]{Definition}
\newtheorem{example}[theorem]{Example}
\newtheorem{remark}[theorem]{Remark}
\newtheorem{conjecture}[theorem]{Conjecture}

\newcommand \bth[1] { \begin{theorem}\label{t#1} }
\newcommand \ble[1] { \begin{lemma}\label{l#1} }

\newcommand \bpr[1] { \begin{proposition}\label{p#1} }
\newcommand \bco[1] { \begin{corollary}\label{c#1} }
\newcommand \bde[1] { \begin{definition}\label{d#1}\rm }
\newcommand \bex[1] { \begin{example}\label{e#1}\rm }
\newcommand \bre[1] { \begin{remark}\label{r#1}\rm }
\newcommand \bcj[1] { \begin{conjecture}\label{j#1}\rm }

\renewcommand {\eth} { \end{theorem} }
\newcommand {\ele} { \end{lemma} }

\newcommand {\epr} { \end{proposition} }
\newcommand {\eco} { \end{corollary} }
\newcommand {\ede} { \end{definition} }
\newcommand {\eex} { \end{example} }
\newcommand {\ere} { \end{remark} }
\newcommand {\ecj} { \end{conjecture} }

\newcommand {\enota} { \end{notation} }
\newcommand \thref[1]{Theorem \ref{t#1}}
\newcommand \leref[1]{Lemma \ref{l#1}}
\newcommand \prref[1]{Proposition \ref{p#1}}
\newcommand \coref[1]{Corollary \ref{c#1}}

\newcommand \deref[1]{Definition \ref{d#1}}
\newcommand \exref[1]{Example \ref{e#1}}
\newcommand \reref[1]{Remark \ref{r#1}}
\newcommand \lb[1]{\label{#1}}
\newtheorem{thmintro}{Theorem}

\def \KK {{\mathbb C}}
\def \KK {{\mathbb K}}
\def \Gset {{\mathbb G}}
\def \Zset {{\mathbb Z}}
\def \Cset {{\mathbb C}}
\def \Nset {{\mathbb N}}
\def \Qset {{\mathbb Q}}

\def \cA {\mathcal{A}}
\def \cC {\mathcal{C}}
\def \cAt {{\mathcal{A}}^{1/2}}
\def \A {{\boldsymbol{\mathsf A}}} 
\def \C {{\boldsymbol{\mathsf C}}} 
\def \R {{\boldsymbol{\mathsf R}}}
\def \U {{\boldsymbol{\mathsf U}}} 

\def \B  {{\widetilde{B}}}
\newcommand \ex {{\mathbf{ex}}}
\newcommand \fr {{\mathbf{fr}}}
\newcommand \inv {{\mathbf{inv}}}
\newcommand \ninv {{\mathbf{ninv}}}

\def \x {{\wt{\mathbf{x}}}} 

\def \SS {{\mathcal{S}}}
\def \LL {{\mathcal{L}}}

\def \TT {{\mathcal{T}}} 

\def \OO {{\mathcal{O}}}

\def \FF {{\mathcal{F}}}

\def \al {\alpha}
\def \be {\beta}
\def \vpi {\varpi}
\def \la {\lambda}

\def \Om {\Omega}
\def \ga {\gamma}

\def \ep {\varepsilon}

\def \mt  {\mapsto}

\def \hra {\hookrightarrow}

\def \rcor {\rangle}
\def \lcor {\langle}

\def \ol {\overline}
\def \wt {\widetilde}
\def \wh {\widehat}

\newcommand{\ubi}{\underline{\mathbf{i}}}

\def \id { {\mathrm{id}} }

\def \sign { {\mathrm{sign}} }

\def \g  {\mathfrak{g}}   

\def \h  {\mathfrak{h}}

\def \n  {\mathfrak{n}}

\def \b  {\mathfrak{b}}

\DeclareMathOperator \Span { {\mathrm{Span}} }
\DeclareMathOperator \Aut { {\mathrm{Aut}} }
\DeclareMathOperator \Oz { {\mathrm{Oz}}}

\DeclareMathOperator \charr { {\mathrm{char}} }

\DeclareMathOperator \diag { {\mathrm{diag}} }
\DeclareMathOperator \Ker { {\mathrm{Ker}} }

\DeclareMathOperator \End { {\mathrm{End}} }

\DeclareMathOperator \Fr { {\mathrm{Fr}}}

\renewcommand \Im { {\mathrm{Im}} }

\renewcommand \max { {\mathrm{max}} }

\newcommand \Fract {{\mathrm{Fract}}}
\begin{document}
\title
[Automorphisms of root of unity quantum cluster algebras]
{Frobenius and ozone automorphism groups \\ of root of unity  quantum cluster algebras}
\author[Y. Wu]{Yunmeng Wu}
\address{Department of Mathematics, Northeastern University, Boston, MA 02115} 
\email{wu.yunm@northeastern.edu}
\address{}
\author[M. Yakimov]{Milen Yakimov}
\address{
Department of Mathematics, Northeastern University, Boston, MA 02115 and International Center for Mathematical Sciences, Institute of Mathematics and Informatics \\
Bulgarian Academy of Sciences \\ 
Acad. G. Bonchev Str., Bl. 8 \\
Sofia 1113, Bulgaria}
\email{m.yakimov@northeastern.edu}
\begin{abstract} Many important families of algebras $R$ in the areas of quantum groups, quantum topology and noncommutative algebra are root of unity quantum algebras, equipped with Frobenius centers $Z_{\Fr}(R) \subseteq Z(R)$ and quantum Frobenius maps $\Fr : R_{\mathrm{class}} \stackrel{\cong}{\longrightarrow} Z_{\Fr}(R)$. We explicitly classify the Frobenius automorphism groups $\Aut(R/Z_{\Fr}(R))$ and the ozone groups $\Oz(R):= \Aut(R/Z(R))$, and prove that all endomorphisms of $R$ fixing its Frobenius center $Z_{\Fr}(R)$ are automorphisms for all algebras $R$ that come from partial root of unity quantum cluster and upper quantum cluster algebras (which have been proved to encompass a vast majority of known quantum algebras). We prove three main theorems that give the above classification results from three different standpoints. As applications, we classify the Frobenius and ozone automorphism groups of all quantum function algebras of simple algebraic groups $G$ not of type $F_4$ or $E_8$ and all quantum unipotent cells at roots of unity for arbitrary symmetrizable Kac--Moody algebras. 
\end{abstract}
\maketitle
\sectionnew{Introduction}
\lb{Intro}
\subsection{Frobenius and ozone automorphism groups of quantum algebras}
\label{1.1}
Throughout the paper $\KK$ will denote a field of arbitrary characteristic. For a $\KK$-algebra $R$ (which is noncommutative in general) and a subalgebra $C \subseteq R$, we consider 
the {\em{relative endomorphism monoid}}
\[
\End(R/C) := \{ \phi \in \End R \mid 
\phi(r) = r, \; \forall r \in C \}
\]
and the {\em{relative automorphism group}}
\[
\Aut(R/C) := \{ \phi \in \Aut R \mid 
\phi(r) = r, \; \forall r \in C \}.
\]
For example, the focal points in Galois theory are the groups $\Aut(F/E)$ for field extensions $F/E$. 

In many areas of mathematics and mathematical physics we have examples of quantum algebras at roots of unity $R$ that are quantizations of commutative algebras $R_{\mathrm{class}}$ and for which there is an injective {\em{quantum Frobenius homomorphism}} from the classical algebra $R_{\mathrm{class}}$ to the center of the quantum one $R$:
\[
R_{\mathrm{class}} \stackrel{\Fr}{\longrightarrow} Z(R) \hookrightarrow R.
\]
Its image 
\[
Z_{\Fr}(R) \cong R_{\mathrm{class}}
\]
is called the {\em{Frobenius center}} of $R$. Quantum Frobenius maps play a key role in the theory of quantum groups, quantum topology, and cluster algebras \cite{BW,DK,DKP,DKP95,DCL94,FG,FKL,GJS,KL,Lus2,NTY2}.

In these situations, it is important to know the symmetries of the quantum object fixing its Frobenius center, i.e., the group of symmetries of the quantum object fixing the classical one embedded in it. This amounts to describing the {\em{Frobenius automorphism group}}
\[
\Aut(R/Z_{\Fr}(R)).
\]
The problem is wide open.

Chan, Gaddis, Won, and Zhang \cite{CGWZ,CGWZ2} introduced the concept of {\em{ozone group}}, which is the automorphism group of a noncommutative algebra $R$ fixing its full center 
\[
\Oz(R):= \Aut(R/Z(R)).
\]
They proved that it plays an important role in the study of Artin--Schelter regular algebras, Calabi--Yau algebras and extensions of the Auslander isomorphism theorem to actions of Hopf algebras. This led to a substantial amount of research \cite{GY,LWZ,Z}. Among many key developments, Liu, Wu and Zhu \cite{LWZ} proved a conjecture of \cite{CGWZ2} that noetherian PI Artin--Schelter regular algebras have (finite) abelian ozone groups. 
There is very little known about computing the groups $\Oz(R)$ in concrete important examples. The cases of skew polynomial algebras, the three-dimensional PI Sklyanin algebras and some quadratic PI Artin--Schelter
regular algebras of global dimension 3 were considered in \cite{CGWZ,CGWZ2}, and certain algebras with 3 generators in \cite{GY}. Many of these groups are trivial. Apart from these examples, ozone groups are unknown. 

Recall that the Dixmier problem asks whether every endomorphism of a Weyl algebra is an automorphism, which is now known to be false in dimensions $\geq 3$. Motivated by it, one can ask if any endomorphism of a noncommutative algebra $R$ that fixes its Frobenius center $Z_{\Fr}(R)$ (or more restrictively, its full center $Z(R)$) is an automorphism. 

In this paper we fully resolve the following three problems:
\begin{enumerate}
\item computation of the Frobenius automorphism group $\Aut(R/Z_{\Fr}(R))$,
\item computation of the ozone group $\Oz(R)= \Aut(R/Z(R))$, and
\item proving that $\End(R/Z_{\Fr}(R)) = \Aut(R/Z_{\Fr}(R))$
\end{enumerate}
for arbitrary root of unity quantum cluster algebras. Multiple descriptions of the automorphism groups are given from different points of view. The results have various applications since many important algebras in the theory of quantum groups and quantum topology, as well as noncommutative algebra, are special cases. 

\subsection{Quantum cluster algebras} 
\label{1.2}
In 2002 Fomin and Zelevinsky introduced the broad axiomatically defined class of cluster algebras \cite{FZ1} for the purposes of studying canonical bases and total positivity. Soon, cluster algebras turned into a major driving force in many areas of mathematics and mathematical physics. Berenstein and Zelevinsky developed the corresponding quantum notion \cite{BerZe}. 
Numerous important families of algebras in various areas of mathematics are special cases of quantum cluster algebras, e.g., \cite{GLS,GY1,OQY26,QY} in the area of quantum groups, \cite{MaQin,M} in quantum topology, and \cite{GY0,GY2} in noncommutative algebra. By a general argument, the specializations of these algebras at roots of unity are root of unity quantum cluster algebras \cite[Theorem 5.11]{NTY2}. 

To state our results, fix a field $\KK$ of characteristic $\charr \KK$ and a primitive $\ell$-th root of unity $\ep^{1/2} \in \KK$ for a positive integer $\ell$. We have a pair of {\em{root of unity quantum and upper quantum cluster algebras}} 
\[
\A_\ep(M_\ep, \B) \subseteq \U_\ep(M_\ep, \B),
\]
associated to a root of unity toric frame $M_\ep$ and an exchange matrix $\B$ whose principal part is skew-symmetrized by a diagonal matrix $D$ with relatively prime positive integral diagonal entries, see Sect. \ref{2.2}.
The first algebra is generated by all cluster variables and the inverses of the frozen ones, while the second algebra is the intersection of the corresponding quantum tori. More generally, having a collection of seeds $\Theta$ which is a nerve in the sense of \cite{F} (see \deref{nerve}) and a subset of indices $\inv$ of the set of frozen ones $\fr$, we consider the {\em{partial root of unity quantum and upper quantum cluster algebras}} 
\begin{equation}
    \label{AU}
\A_\ep(M_\ep, \B, \inv,\Theta) \subseteq \U_\ep(M_\ep, \B, \inv,\Theta),
\end{equation}
obtained from the seeds in $\Theta$ and where only the frozen variables in the subset $\inv$ of all frozen variables $\fr$ are inverted, while the other frozen variables are not. The introduction of a subset of seeds $\Theta$ adds a lot of flexibility in concrete applications, since it is a lot easier to prove that a given algebra is isomorphic to one of the above type for a finite set $\Theta$. The root of unity quantum cluster algebras in \eqref{AU} have canonical central subalgebras (Frobenius centers)
\[
\C\A_\ep(M_\ep, \B, \inv,\Theta) \subseteq \C\U_\ep(M_\ep, \B, \inv,\Theta),
\]
obtained in a similar way to the ambient algebras but by using the $\ell$-th powers of the cluster/frozen variables instead, see Sect. \ref{2.2}. If $\ell$ is a positive integer satisfying the condition
\smallskip

\begin{description}
\item[{\bf{(Coprime)}}] {\em{$\ell$ is odd and coprime to $\charr \KK$ and the diagonal entries of the skew-symmetrizing matrix $D$ for the principal part of the exchange matrix $\B$,}}
\end{description}
\smallskip
\noindent
then there are quantum Frobenius isomorphisms \cite{HLY,NTY2} 
\begin{align*}
\A(\B, \inv, \Theta) &\stackrel{\simeq}{\longrightarrow} \C\A_\ep(M_\ep, \B, \inv, \Theta), 
\\
\U(\B, \inv, \Theta) &\stackrel{\simeq}{\longrightarrow} \C\U_\ep(M_\ep, \B, \inv, \Theta),
\end{align*}
see \thref{cent-isom}(ii)--(iii).
(Coprime) is a condition on a root of unity quantum
cluster algebra and not on individual seeds \cite[Remark 4.3]{NTY2}, meaning that, if it is satisfied for one seed, then it is satisfied for all seeds obtained from it by mutation. 

Throughout the rest of the introduction we denote 
\[
\R_\ep(M_\ep, \B, \inv, \Theta):= \A_\ep(M_\ep, \B, \inv, \Theta) \quad \mbox{or}
\quad \U_\ep(M_\ep, \B, \inv, \Theta),
\]
and, respectively,
\[
\C\R_\ep(M_\ep, \B, \inv, \Theta):= \C\A_\ep(M_\ep, \B, \inv, \Theta) \quad \mbox{or}
\quad \C\U_\ep(M_\ep, \B, \inv, \Theta).
\]
All theorems will hold in both cases. 

For an integer matrix $C$, we will denote by $\ol{C}$ its reduction modulo $\ell$. 
There is a faithful action of the kernel $\Ker \big( \, \ol{\B}^\top \big)$ in $(\Zset/\ell)^N$ on $\R_\ep(M_\ep,\B, \inv, \Theta)$ given by 
\[
\theta \cdot M_\ep(f) := \ep^{\theta \cdot f} M_\ep(f), \quad \forall 
\theta \in \Ker \big( \, \ol{\B}^\top \big), f \in \Zset^N,
\]
see \prref{ep-act}. Here and below $N$ denotes the number of rows of $\B$. Our first classification results work under weaker 
assumptions than those in condition (Coprime). 

\begin{thmintro} 
\label{thmA}
 Assume that $\ell$ is an odd positive integer coprime to $\charr \KK$, $\ep^{1/2} \in \KK$ is a primitive $\ell$-th root of unity, $\Theta$ is a connected collection of seeds of the root of unity quantum cluster algebra $\A_\ep(M_\ep, \B)$ that is a nerve, and $\inv \subseteq \fr$ is arbitrary. Then the following hold:
 \begin{enumerate}
 \item[(i)] The Frobenius endomorphism monoids and Frobenius automorphism groups of the corresponding partial root of unity quantum and upper quantum cluster algebras are given by 
\begin{align*}
&\End \big( \R_\ep(M_\ep, \B, \inv, \Theta)/
\C\R_\ep(M_\ep, \B, \inv, \Theta) \big) = \\
& \Aut \big( \R_\ep(M_\ep, \B, \inv, \Theta)/
\C\R_\ep(M_\ep, \B, \inv, \Theta) \big) = \Ker \big( \, \ol{\B}^\top \big),
\end{align*}
where the second equality means that every automorphism is realized
as the action $(\theta \cdot)$ of an element of the corresponding group in the right-hand side.
\item[(ii)] The ozone groups of these algebras are given by 
\[
\Oz \big( \R_\ep(M_\ep, \B, \inv, \Theta) \big) = 
\Ker \big( \, \ol{\B}^\top \big) \cap \big( \ol{\Ker \Om_{M_\ep}} \big)^\perp,
\]
where $\Om_{M_\ep}$ is the $N \times N$ integer matrix of the frame $M_\ep$ and the orthogonal complement is computed with respect to the dot product on $(\Zset/\ell)^N$.
\end{enumerate}
\end{thmintro}

A more general result for the automorphism groups of algebras sandwiched between  $\A_\ep(M_\ep, \B, \inv, \Theta)$
and $ \U_\ep(M_\ep, \B, \inv, \Theta)$ is obtained in 
\thref{aut-qcl}. The theorem classifies the Frobenius and ozone automorphism groups of all 
partial root of unity quantum and upper quantum cluster algebras in terms of the initial exchange matrix $\B$ and fully resolves the three problems in Sect. \ref{1.1}.

In concrete cases when an algebra is shown to be isomorphic to an algebra of the form $\R_\ep(M_\ep, \B, \inv, \Theta)$, computations with the exchange matrix can be complicated because of its large size. Our next main theorem provides a simpler answer to the classification problem from this practical perspective when condition (Coprime) holds.

\begin{thmintro} 
\label{thmB}
Assume that $\ell$ is a positive integer satisfying condition (Coprime), $\ep^{1/2} \in \KK$ is a primitive $\ell$-th root of unity, $\Theta$ is a connected collection of seeds of the root of unity quantum cluster algebra $\A_\ep(M_\ep, \B)$ that is a nerve, and $\inv \subseteq \fr$ is arbitrary. 
\begin{enumerate}
\item[(i)] Assume that $T$ is a split torus over $\KK$ of rank 
$|\fr|$ acting on $\R_\ep(M_\ep, \B, \inv, \Theta)$
such that the $\ell$-torsion subgroup $T(\KK)[\ell]$ of the group 
of $\KK$-points $T(\KK) \cong (\KK^\times)^{|\fr|}$ acts faithfully and the cluster variables $M'_\ep(e_i)$ are 
$T(\KK)$-eigenvectors for all $(M'_\ep, \B') \in \Theta$ and $i \in [1,N]$.

The Frobenius automorphism groups of the corresponding partial root of unity quantum and upper quantum cluster algebras are given by  
\[
\Aut \big(
\R_\ep(M_\ep, \B, \inv, \Theta)/
\C\R_\ep(M_\ep, \B, \inv, \Theta) \big) = 
T(\KK)[\ell] = C_\ell^{|\fr|} \cong (\Zset/\ell)^{|\fr|}.
\]
Here $C_\ell \subset \KK^\times$ is the group of $\ell$-th roots of unity and 
the equalities mean that every element of the relative automorphism 
groups comes from the action of the $\ell$-torsion subgroup of $T$.

The ozone groups of these algebras are given by 
\[
\Oz \big( \R_\ep(M_\ep, \B, \inv, \Theta) \big) =
\{ t \in C_\ell^{|\fr|} \mid t \cdot M_\ep(f) = M_\ep(f), \; \forall 
f \in \Ker (\Om_{M_\ep})\}.
\]
\item[(ii)] There is a canonical split torus over $\KK$ of rank $|\fr|$ with the properties in part (i) with an action on $\R_\ep(M_\ep, \B, \inv, \Theta)$ 
explicitly constructed in \prref{Caction}. In particular, 
\[
\Ker \big( \, \ol{\B}^\top \big) \cong (\Zset/ \ell)^{|\fr|}.
\]
\end{enumerate}
\end{thmintro}

The first part of the theorem is extremely useful in applications. If a given algebra is proved to be isomorphic to a partial root of unity quantum or upper quantum cluster algebra and condition (Coprime) is satisfied, then we do not need to know the exact structure of the exchange matrix $\B$ to classify its Frobenius and ozone automorphism groups. We can just use a canonical torus action on the algebra and check that the cluster variables for the seeds in the set $\Theta$ are eigenvectors for the group of $\KK$-points of the torus. The second part of the theorem shows that the torus action in it is not artificial; it always exists if condition (Coprime) holds.  

In \exref{G2-noncoprime} we show that the Frobenius automorphism groups of the algebras $\R_\ep(M_\ep, \B, \inv, \Theta)$ need not be isomorphic to $(\Zset/ \ell)^{|\fr|}$ if condition (Coprime) is violated. This illustrates a situation 
when Theorem \ref{thmA} is applicable but Theorem \ref{thmB} is not. 

Chan, Gaddis, Won, and Zhang \cite{CGWZ2} proved that for every PI ring which is a domain, its ozone group consists of the automorphisms 
given by conjugations by its normal elements. In general, this does not classify the group on its own since PI rings have a huge number 
of normal elements and it is very hard to get any classification results for the multiplicative sets of their normal elements except for very simple examples.  
Our third main result proves that in the case of partial root of
unity quantum and upper quantum cluster algebras one can make this relation very explicit in terms of the frozen variables. 

\begin{thmintro} 
\label{thmC}
Assume that $\KK$, $\B$, $\ell$, $\ep^{1/2}$, $\Theta$, $\inv$, and $\fr$ are as in Theorem \ref{thmB}. The ozone automorphism groups 
of the partial root of unity quantum and upper quantum cluster algebras $\R_\ep(M_\ep, \B, \inv, \Theta)$ are
\[
\Oz \big( \R_\ep(M_\ep, \B, \inv, \Theta) \big) \cong \Om_{M_\ep} \Big( \bigoplus_{i \in \fr} (\Zset/\ell) \ol{e_i} \Big) \subseteq (\Zset/\ell)^N,
\]
where $\{ \ol{e_i} \}_{i=1}^N$ is the standard basis of $(\Zset/\ell)^N$. 

All elements of these groups are realized as conjugations by monomials in the frozen variables 
\[
M_\ep \Big( \sum_{i \in \fr} n_i e_i \Big), \quad \mbox{for} \quad n_i \in [0, \ell -1].
\]
\end{thmintro} 

The proofs of the three main theorems are dependent: Theorems \ref{thmB} and \ref{thmC} use Theorem \ref{thmA}.

\subsection{Application} The three main results from the previous subsection have broad applications to the classification 
of the Frobenius and ozone automorphism groups of important families of quantum algebras at roots of unity. 

Let $G$ be a connected, simply connected, complex simple algebraic group and $\ell$ a positive odd integer that is coprime to the 
symmetrizing integers for its Cartan matrix. De Concini and Lyubashenko \cite{DCL94} 
constructed a (big) quantized coordinate ring of $G$ 
at a primitive $\ell$-th root of unity $\ep$ over $\Cset$ and a Frobenius isomorphism 
\[
\Fr_\ep : \OO(G) \stackrel{\cong}{\longrightarrow} \cC_\ep(G) \subset Z( \OO_\ep(G)),
\]
see Sect. \ref{qgroups} for details. 

\begin{thmintro} 
\label{thmD}
Let $G$ be a connected, simply connected, complex simple algebraic group of rank $r$ not of type $F_4$ or $E_8$, $\ell$ be an odd positive integer not divisible by 3 if $G$ is of type $G_2$, and $\ep^{1/2}$ be a primitive $\ell$-th root of unity. 

\begin{enumerate}
\item[(i)] The Frobenius automorphism group of $\OO_\ep(G)$  is given by 
\[
\Aut(\OO_\ep(G)/\cC_\ep(G)) \cong (\Zset/\ell)^{2r}.
\]

Every endomorphism of $\OO_\ep(G)$ fixing $\cC_\ep(G)$ is an automorphism.

\item[(ii)] The ozone group of $\OO_\ep(G)$ is isomorphic to the diagonal subgroup of the group in the first part:
\[
\Oz(\OO_\ep(G)) \cong (\Zset/\ell)^r.
\]
The elements of the group $\Oz(\OO_\ep(G))$ are precisely the automorphisms of $\OO_\ep(G)$ of the form 
\[
r \mt 
\Delta_{w_0 \la, \la} r \Delta_{w_0 \la, \la}^{-1} \quad \mbox{for} \quad\la = n_1 \varpi_1 + \cdots + n_r \varpi_r, \; \mbox{with} \; n_i \in [0, \ell -1],
\]
where $\Delta_{w_0 \la, \la}$ are the standard extremal matrix coefficients given by \eqref{la}.
\end{enumerate}
\end{thmintro}

The restriction to types different from $F_4$ and $E_8$ is needed because of the results in \cite{OQY26} on integral forms of quantum cluster algebras. 

Let us consider, more generally, a symmetrizable Kac--Moody algebra $\g$ and an odd positive integer $\ell$ 
that is coprime to the symmetrizing integers for its Cartan matrix. To every Weyl group element $w \in W$, 
following De Concini--Kac--Procesi \cite{DKP95}, Lusztig \cite{Lus93} and Gei\ss--Leclerc--Schr\"oer \cite{GLS}, one associates the 
{\em{quantum unipotent cell}} (also called {\em{quantum Schubert cell}}) $A_\ep(\n_+(w))$ and its Frobenius center
\[
C_\ep(\n_+(w)) \subset Z(A_\ep(\n_+(w))), 
\]
see Sect. \ref{unip-Schubert} for details. 
Our last main result classifies the Frobenius and ozone automorphism groups of all quantum unipotent cells at roots of unity.

\begin{thmintro} 
\label{thmE}
Let $\g$ be a symmetrizable Kac--Moody algebra, $w \in W$, and $\ep$ a primitive $\ell$-th root of unity, where $\ell$ is odd and coprime to the symmetrizing integers for the Cartan matrix of $\g$. 
\begin{enumerate}
    \item[(i)] The Frobenius automorphism group of the quantum unipotent cell $A_\ep(\n_+(w))$ is  
\[
\Aut(A_\ep(\n_+(w))/C_\ep(\n_+(w))) \cong (\Zset/\ell)^{|\SS(w)|}, 
\]
where $\SS(w)$ denotes the support of $w$, i.e., the collection of indices of the simple reflections that appear in one (and thus in any) reduced expression of $w$. 

Every endomorphism of $A_\ep(\n_+(w))$ fixing $C_\ep(\n_+(w))$ is an automorphism.
    \item[(ii)] The ozone group of $A_\ep(\n_+(w))$ is
\[
\Oz(A_\ep(\n_+(w))) \cong \wh{Q_{\SS(w), \ell}},
\]
where the right-hand side denotes the character group of the finite abelian group
\[
Q_{\SS(w),\ell}:=
(1+w^{-1})Q_{\SS(w)}/\big( (1+w^{-1})Q_{\SS(w)} \cap \ell Q_{\SS(w)} \big)
\]
and $Q_{\SS(w)}$ is the sublattice of the root lattice generated by the simple roots indexed by $\SS(w)$. 

The elements of the group $\Oz(A_\ep(\n_+(w)))$ are conjugations by monomials 
in the unipotent quantum minors $D_{\vpi_i, w \vpi_i}$ for $i \in \SS(w)$, 
given by \eqref{unip-qcomm}.
\end{enumerate}
\end{thmintro}

In \thref{complex} we use Theorem \ref{thmE} to illustrate the very interesting case of the ozone group of the half-quantum group $U_\ep(\n_+)$ for all finite dimensional complex simple Lie algebras $\g$.  

\subsection{Notations}
The character lattice of a split torus $T \cong \Gset_m^r$ will be denoted by $X(T)$. The character group of a 
finite group $G$ will be denoted by $\wh{G}$. 

For a positive integer $N$, we will denote $[1,N]:= \{1, \ldots, N \}$. Set $\Nset := \{0, 1, \ldots\}$.
The standard basis of $\Zset^N$ will be denoted by $\{e_1, \ldots, e_N\}$. Given $\nu = (\nu_1, \ldots, \nu_N) \in \Zset^N$ or $(\Zset/\ell)^N$, set
\begin{equation}
\label{j-coord}
(\nu)_j:= \nu_j, \quad \forall j \in [1,N].
\end{equation}
Let $ |\nu|, [\nu]_\pm \in \Zset^N$ be the unique vectors such that 
\[
(|\nu|)_j := |\nu_j|, 
\quad 
([\nu]_\pm)_j := \max (\pm \nu_j, 0), \; \; \forall j \in [1,N].
\]
We will denote by 
\begin{equation}
\label{dot-prod}
\nu \cdot \pi := (\nu)_1 (\pi)_1 + \cdots + (\nu)_N (\pi)_N
\end{equation}
the dot product on $\Zset^N$ and $(\Zset/\ell)^N$. 

For two finite sets of integers $I$ and $J$ and a commutative ring $C$, we will denote by $M_{I \times J}(C)$ the rectangular matrices with entries in $C$ whose rows are indexed by $I$ and columns by $J$. Set $M_S(C):= M_{S \times S}(C)$. If $S =[1,N]$, we will abbreviate $S$ in the matrix notation by $N$.

\subsection*{Terminology} For a map $\phi : X \to X$ and a subset $Y \subseteq X$, we will say that $\phi$ preserves $Y$ if $\phi(Y) \subseteq Y$ and that $\phi$ fixes $Y$ if $\phi(y) =y$ for all $y \in Y$. 

\subsection*{AI statement} No AI tools have been used in any stage of the project, including research and editing.  

\subsection*{Acknowledgments} 
The research of Y.W. was partially supported by NSF grants DMS-2601486
and DMS-1645877. The research of M.Y. was supported by NSF grant DMS-2601486, Bulgarian Science Fund grant KP-06-N92/5, the Ministry of Education and Science grant DO1-239/10.12.2024, and Simons Foundation grant SFI-MPS-T-Institutes-00007697.
\sectionnew{Preliminaries on cluster algebras}
\lb{prelim}
\subsection{Cluster algebras of geometric type}
\label{2.1}
Here we review the basics of cluster algebras of geometric type following Berenstein, Fomin and Zelevinsky \cite{BFZ05,FZ1}.
Fix a positive integer $N$. Let $\FF$ be a purely transcendental field extension of $\KK$ 
of transcendence degree $N$. Fix a set of {\em{exchangeable indices}}  $\ex\subseteq [1,N]$
and denote the {\em{frozen indices}}
\[
\fr := [1,N] \setminus \ex.
\]
We say that a pair $(\x, \B)$ is a {\em{seed}} if 
\begin{enumerate}
\item[(i)] $\x = \{x_1, \dots, x_N\}$ is a transcendence basis of $\FF$ over $\KK$;
\item[(ii)] $\B = (b_{ij})
\in M_{N\times \ex}(\Zset)$ and its principal $\ex\times \ex$ submatrix $B$
is skew-symmetrizable; that is, there exists a diagonal matrix $D=\diag(d_j; j \in \ex)$ with positive relatively prime entries $d_j \in \Zset_+$ such that $(D B)^\top = -D B$.
\end{enumerate}
\bre{cv} We will call $\{x_1, \ldots, x_N\}$ {\em{cluster variables}}, while $\{x_k \mid k \in \ex\}$ will be called {\em{exchangeable}} ({\em{mutable}}) ones and $\{x_k \mid k \in \fr\}$ frozen ones. This differs from some conventions where the second set is called the set of cluster variables and the third one the set of frozen variables or coefficients. 
\ere


The mutation of the {\em{exchange matrix}} $\B$ in direction $k \in \ex$ is the matrix 
\begin{equation}
\label{mutation}
\mu_k(\B) = (b_{ij}'):=
\begin{cases}
-b_{ij}, & \text{if } i = k, \text{ or } j=k \\
b_{ij} + \frac{|b_{ik}|b_{kj} + b_{ik}|b_{kj}|}{2}, & \text{otherwise.}
\end{cases}
\end{equation}
For a seed $(\x, \B)$, its mutations are the seeds $\mu_k (\x, \B) := (\x', \mu_k(\B))$, where 
\[
\x' := \{x_k'\} \cup \x \setminus \{x_k\} 
\quad
\mbox{and}
\quad
x_k x_k' := \prod_{b_{ik}>0}x_i^{b_{ik}} + \prod_{b_{ik}<0}x_i^{-b_{ik}}
\]
for $k \in \ex$. Given $s = \pm$, define the matrices $E_s \in M_N(\Zset)$ and $F_s \in M_\ex(\Zset)$ to be
\[
E_s := (e_{ij}) =
\begin{cases}
\delta_{ij} &  \text{if } j \neq k \\
-1 & \text{if } i = j = k \\
\max(0,-sb_{ik}) & \text{if } i \neq j = k,
\end{cases} \; \; 
F_s := (f_{ij})=
\begin{cases}
\delta_{ij} &  \text{ if } i \neq k \\
-1 & \text{ if } i = j = k \\
\max(0,sb_{kj}) & \text{ if } j \neq i = k.
\end{cases}
\]
We have $\mu_k(\B)= E_s \B F_s$ for both $s = \pm$.

Mutation is an involution. Two seeds are {\em{mutation equivalent}}, $(\x', \B')\sim (\x'', \B'')$, if one can be obtained from the other by a finite sequence of mutations.
Any seed which is mutation equivalent to $(\x, \B)$ contains the {\em{frozen variables}}, $x_i$, $i \in \fr$. 

The {\em{exchange graph of the cluster algebra}} is the labeled graph with vertex set the set of distinct seeds and edges representing mutation, labeled by the elements of $\ex$.

Fix a decomposition of the set of frozen indices: 
\[
\fr = \inv \sqcup \ninv,
\]
where the frozen variables $x_i$ for $i \in \inv$ will be inverted and those $x_i$ for $i \in \ninv$ will not. Given a seed $(\x',\B')\sim (\x,\B)$, denote the {\em{Laurent polynomial algebra}}
\[
\LL(\x') = \KK[(x'_k)^{\pm 1}; 1 \leq k \leq N]
\]
and the {\em{mixed torus-affine subalgebra}}
\[
\LL(\x')_\geq = \KK[(x'_k)^{\pm 1}, x'_i; k \in \ex \sqcup \inv, i \in \ninv].
\]
\bde{CAs} \cite{BFZ05,FZ1}
\hfill 
\begin{enumerate}
\item[(i)] The {\em{partially compactified cluster algebra}} $\A(\B, \inv)$ is the $\KK$-subalgebra of $\FF$ generated by 
the cluster variables in all the seeds $(\x',\B')\sim( \x, \B)$ together with $\{ x_i^{-1} \mid i \in \inv\}$.
\item[(ii)] The {\em{partially compactified upper cluster algebra}} $\U(\B, \inv)$ is the intersection 
\[
\U(\B,\inv):= \bigcap_{(\x', \B') \sim (\x, \B)} \LL(\x')_\geq.
\]
\end{enumerate}
\ede

The {\em{Laurent phenomenon}} \cite{FZlp} gives that $\A(\B, \inv) \subseteq \U(\B, \inv)$. 

The classical cluster algebra and upper cluster algebra correspond to the case when $\ninv = \varnothing$, in which case they are denoted by 
$\A(\B)$ and $\U(\B)$, respectively. They are obtained by localization from their partially compactified counterparts 
\[
\A(\B) = \A(\B, \inv)[x_i^{-1}; i \in \ninv], \quad 
\U(\B) = \U(\B, \inv)[x_i^{-1}; i \in \ninv].
\]

More generally, we associate two algebras to a {\em{connected collection of seeds}} $\Theta$ of the cluster algebra \cite{F}. Here and below, by connected we mean connected in the exchange graph, i.e., every two seeds in $\Theta$ can be reached from one another by mutations. 
The algebra  
$\A(\B,\inv, \Theta)$ is the subalgebra of the cluster algebra $\A(\B,\inv)$ generated by the cluster variables in the seeds in $\Theta$ and the inverses of the frozen variables in $\inv$. We also consider the partial intersection 
\[
\U(\B,\inv, \Theta):= \bigcap_{(\x', \B') \in \Theta} \LL(\x')_\geq
\] see \cite{F,HLY}. The Laurent phenomenon implies that for all collections $\Theta$, we have 
\[
\A(\B,\inv, \Theta) \subseteq \A(\B,\inv) \subseteq \U(\B,\inv) \subseteq \U(\B,\inv, \Theta). 
\]
\subsection{Root of unity quantum cluster algebras}
\label{2.2}
Here we review the basics of root of unity quantum cluster algebras following \cite{NTY2}.

Fix a positive integer $\ell$, coprime to $\charr \KK$, and a primitive $\ell$-th root of unity $\ep^{1/2} \in \KK$. By abuse of notation, for every skew-symmetric bilinear form $\Omega : \Zset^N \times \Zset^N \to \Zset/\ell$, we will denote by the same symbol its matrix 
\[
\Omega:=(\Omega(e_i,e_j))_{i,j=1}^N \in M_N(\Zset/\ell). 
\]
Such a bilinear form gives rise to the {\em root of unity based quantum torus}
\begin{equation}
\label{Tep}
\TT_\ep(\Omega) := \Span_{\KK} \{ x^f \mid f \in \Zset^N \}, \; \; 
\mbox{where} \; \; x^f x^g = \ep^{\Omega(f,g)/2} x^{f+g}, \; \; \forall f, g \in \Zset^N
\end{equation}
and the {\em{root of unity mixed quantum torus-affine algebra}}
\[
\TT_\ep(\Om)_{\geq} :=\KK \lcor x_k^{\pm 1}, \, x_i  :  k \in \ex \sqcup \inv,  i \in \ninv \rcor \subset \TT_\ep(\Om).
\]

A {\em{root of unity toric frame}} $M_\ep$ of a division algebra $\FF$ over $\KK$ is a map 
\[
M_\ep : \Zset^N \to \FF,
\]
which extends to a $\KK$-algebra embedding $\phi: \TT_\ep(\Omega) \hra \FF$ such that
\begin{enumerate}
    \item $\phi(x^\nu)=M_\ep(\nu)$ for all $\nu \in \Zset^N$ and 
    \item $\FF \simeq \Fract \left( \TT_\ep(\Omega) \right)$.
\end{enumerate}
The bilinear form is uniquely reconstructed from the toric frame. It is called {\em{the bilinear form of the toric frame}} and denoted by $\Omega_{M_\ep}$. Set 
\[
\TT_\ep(M_\ep)_{\geq} := \phi (\TT_\ep(\Om)_{\geq}).
\]

\bre{cv2} Similarly to \reref{cv}, we will call $\{M_\ep(e_k) \mid k \in [1,N] \}$ {\em{cluster variables}}, while $\{M_\ep(e_k) \mid k \in \ex\}$ will be called {\em{exchangeable}} ({\em{mutable}}) ones and $\{M_\ep(e_k) \mid k \in \fr\}$ frozen ones. Once again, this slightly differs from the conventions of some authors who call the second set the set of cluster variables, and the third one the set of frozen variables or coefficients. 
\ere

We say that a pair $(\Omega, \B)$, consisting of a skew-symmetric bilinear form  $\Omega  : \Zset^N \times \Zset^N \to \Zset/ \ell$ and an exchange matrix $\B \in M_{N \times \ex}(\Zset)$,
is {\em{$\ell$-compatible}} if there exists a diagonal matrix $D := \diag(d_j, j \in \ex)$ with $d_j \in \Zset_+$ 
such that the principal part $B$ of $\B$ is skew-symmetrized by $D$ and 
\begin{equation}
\label{compatible}
\ol{\B}^\top \Omega  = 
\begin{bmatrix}
\, \ol{D} \; 0
\end{bmatrix} \hspace{-0.1cm},
\end{equation}
where, for an integer matrix $C$, $\ol{C}$ denotes its reduction modulo $\ell$.
The \emph{mutation} in direction $k \in \ex$ of an $\ell$-compatible pair is defined to be $\mu_k ( \Omega, \B) := ( \ol{E}_s^\top \Omega \ol{E}_s, E_s \B F_s)$ for $s= \pm$, which is independent of the choice of sign $s$. The mutated pairs are $\ell$-compatible with respect to the same diagonal matrix $D$. 

A {\em root of unity quantum seed} is a pair $(M_\ep, \B)$ consisting of a toric frame and an exchange matrix such that $(\Omega_{M_\ep},\B)$ is an $\ell$-compatible pair. 
Its mutation in direction $k\in\ex$ is the seed
$\mu_k(M_\ep, \B) := (\rho^{ M_\ep}_{b^k, s} M_\ep E_s,  E_s \B F_s)$, where $\rho_{b^k,s}=\rho^{M_\ep}_{b^k,s}$ is the unique automorphism of $\FF$ such that 
\[ 
\rho_{b^k, s}^{M_\ep}(M_\ep(e_j)) = 
\begin{cases}
M_\ep(e_k) + M_\ep(e_k + sb^k) & \text{if } j=k \\
M_\ep(e_j)                & \text{if } j \neq k.
\end{cases}
\]
Here and below, we will denote by $b^k$ the $k$-th column of $\B$.
Mutation of root of unity quantum seeds does not depend on the choice of sign and is an involution. 
Two seeds are mutation equivalent $(M'_\ep, \B') \sim (M_\ep, \B)$ if one can be obtained from the other by a sequence of mutations. The frozen variables are unchanged: 
$M'_\ep(e_i) = M_\ep(e_i)$ for $i \in \fr$.

The {\em{exchange graph of the root of unity quantum cluster algebra}} is the labeled graph with vertex set the set of distinct quantum seeds and edges representing mutation, labeled by the elements of $\ex$.

\bde{rootunitqCA}
\hfill
\begin{enumerate}
\item[(i)]
The {\em partially compactified root of unity quantum cluster algebra} $\A_\ep(M_\ep, \B, \inv)$ is the $\KK$-subalgebra of $\FF$ 
generated by all cluster variables of the seeds $(M'_\ep, \B') \sim (M_\ep, \B)$ and by 
$M_\ep(e_i)^{-1}$ for $i \in \inv$.
\item[(ii)] The 
\emph{partially compactified root of unity upper quantum cluster algebra} is the intersection 
\[
\U_\ep(M_\ep, \B, \inv) := \bigcap_{ (M'_\ep, \B') \sim (M_\ep, \B)} \TT_\ep(M_\ep')_{\geq}.
\]
\end{enumerate}
\ede

Similarly to the classical situation, the root of unity quantum cluster algebras and upper quantum cluster algebras are the ones obtained in the case $\ninv = \varnothing$ and are localizations of the general partially compactified versions:
\begin{align*}
\A_\ep(M_\ep,\B) &=
\A_\ep(M_\ep, \B,\inv)[M_\ep(e_i)^{-1}; i \in \ninv],
\\
\U_\ep(M_\ep,\B) &=
\U_\ep(M_\ep,\B, \inv)[M_\ep(e_i)^{-1}; i \in \ninv].
\end{align*}

Let $\Theta$ be a connected collection of seeds of the root of unity quantum cluster algebra,
where by connected we mean connected in the exchange graph of the root of unity quantum cluster algebra. We associate two algebras to it:
\begin{enumerate}
\item the subalgebra $\A_\ep(M_\ep,\B,\inv, \Theta)$  of the root of unity quantum cluster algebra $\A_\ep(M_\ep,\B,\inv)$ generated by the cluster variables in the seeds in $\Theta$ and the inverses of the frozen variables in $\inv$, and 
\item the overalgebra of the root of unity upper quantum cluster algebra $\U_\ep(M_\ep,\B,\inv)$ given by the intersection
\[
\U_\ep(M_\ep, \B, \inv, \Theta) := \bigcap_{ (M'_\ep, \B') \in \Theta}\TT_\ep(M_\ep')_{\geq},
\]
see \cite[Sect. 3.2]{HLY}. 
\end{enumerate}
They will be called partial root of unity quantum and upper quantum cluster algebras.

By a root of unity quantum analog of the Laurent phenomenon \cite[Theorem 3.10]{NTY2}, $\A_\ep(M_\ep, \B, \inv)\subseteq \U_\ep(M_\ep,  \B, \inv)$. Therefore, we have the chain of embeddings
\[
\A_\ep(M_\ep,\B,\inv, \Theta) \subseteq \A_\ep(M_\ep,\B,\inv) \subseteq \U_\ep(M_\ep,\B,\inv) \subseteq \U_\ep(M_\ep,\B,\inv, \Theta). 
\]

For every seed $(M'_\ep, \B') \sim (M_\ep, \B)$ and $f \in \Zset^N$,
\[
M'_\ep(f)^\ell = 
M'_\ep(\ell f) \in Z(\FF). 
\] 
Thus, we have the (central) mixed torus-affine subalgebra
\[
\TT_\ep( M'_\ep)_{\geq}^\ell:= \KK [ M'_\ep(\ell e_k)^{\pm 1}, \ M'_\ep(\ell e_i)   : 
k \in \ex \sqcup \inv,  i \in \ninv ] \subset Z( \TT_\ep(M'_\ep)_\geq)
\subseteq Z(\FF).
\]
For a connected collection $\Theta$ of seeds of the root of unity quantum cluster algebra, we obtain 
\begin{enumerate}
\item 
the central subalgebra
\[
\C\A_\ep(M_\ep, \B, \inv, \Theta)
\subset Z (\A_\ep(M_\ep, \B, \inv, \Theta))
\subset Z(\FF)
\]
generated by $M'_\ep(\ell e_k)$ for seeds 
$(M'_\ep, \B') \in \Theta$ with $k \in [1,N]$, and by $M_\ep(\ell e_i)^{-1}$ for $i \in \inv$ 
and 
\item 
the central subalgebra  
\[
\C\U_\ep(M_\ep, \B, \inv, \Theta) := \bigcap_{ (M'_\ep, \B') \in \Theta} \TT_\ep(M_\ep')_{\geq}^\ell
\subset Z (\U_\ep(M_\ep, \B, \inv, \Theta)) 
\subset Z(\FF).
\]
\end{enumerate}
These central subalgebras will be called the {\em{Frobenius centers}} of $\A_\ep(M_\ep, \B, \inv, \Theta)$ and 
$\U_\ep(M_\ep, \B, \inv, \Theta)$, respectively. 

For a bilinear form $\Om : \Zset^N \times \Zset^N \to \Zset/\ell$, denote 
\[
\Ker \Om := \{ f \in \Zset^N \mid \Om(f,g) =0, \; \forall g \in \Zset^N \}.
\]
In particular, $\Ker \Om \supseteq (\ell \Zset)^N$, so the index of $\Ker \Om$ 
in $\Zset^N$ is finite. It is straightforward to verify that
\begin{equation}
\label{ZT}
Z(\TT_\ep(M'_\ep)_\geq) = \Span_\KK 
\{ M'_\ep(f) \mid f \in \Ker \Om_{M'_\ep}, (f)_i \geq 0, \; \forall i \in \ninv \}
\supseteq \TT_\ep(M'_\ep)_{\geq}^\ell,
\end{equation}
see \cite[Eq. (6.1)]{HLY} and recall the notation \eqref{j-coord}. 

Recall the condition (Coprime) on the positive integer $\ell$ from the introduction.
We have the following facts for quantum Frobenius maps:

\bth{cent-isom} Assume that $\ell$ is a positive integer satisfying condition (Coprime) and $\ep^{1/2} \in \KK$ is a primitive $\ell$-th root of unity. 
Then the following hold:
\begin{enumerate}
\item[(i)] \cite[Theorem 4.8]{NTY2}
There exists a unique isomorphism of labeled graphs between the exchange graphs of $\A(\B)$ and $\A_\ep(M_\ep,\B)$ sending the seed $(\x, \B)$ to $(M_\ep, \B)$. We will use this isomorphism to identify (connected) collections of seeds of the classical and root of unity quantum cluster algebras.
\item[(ii)] 
There is an isomorphism of $\KK$-algebras (quantum Frobenius map)
\[
\A(\B, \inv, \Theta) \stackrel{\simeq}{\longrightarrow} \C\A_\ep(M_\ep, \B, \inv, \Theta), 
\]
given by $x'_k \mt M'_\ep(\ell e_k)$ for all cluster variables of the seeds in $\Theta$ and 
$x_i^{-1} \mt M_\ep(\ell e_i)^{-1}$ for all inverted frozen variables, $i \in \inv$.
\item[(iii)] \cite[Proposition 3.2]{HLY}
For every connected collection of seeds $\Theta$ of $\A(\B)$, the isomorphism from part (ii) extends to a $\KK$-algebra isomorphism 
(also called quantum Frobenius map)
\[
\U(\B, \inv,\Theta) \stackrel{\simeq}{\longrightarrow} \C\U_\ep(M_\ep, \B, \inv, \Theta).
\]

\item[(iv)] \cite[Lemma 3.3]{HLY} If $R$ is a $\KK$-algebra sandwiched as 
\[
\A_\ep(M_\ep, \B, \inv, \Theta) \subseteq R \subseteq 
\U_\ep(M_\ep, \B, \inv,\Theta),
\]
then the center of $R$ is given by  
\[
Z(R) = R \cap Z(\TT_\ep(M'_\ep)_{\geq}),
\]
for every seed $(M'_\ep, \B')$ in $\Theta$, 
where $Z(\TT_\ep(M'_\ep)_\geq)$ is described in \eqref{ZT}.

\item[(v)] \cite[Theorem B]{HLY} If 
$\U_\ep(M_\ep, \B, \inv,\Theta)$ is a finitely generated $\KK$-algebra, then $\U_\ep(M_\ep, \B, \inv,\Theta)$ is a finitely generated 
$\C\U_\ep(M_\ep, \B, \inv,\Theta)$-module and the classical partial upper cluster algebra $\U(\B, \inv,\Theta)$ is a finitely generated $\KK$-algebra. 
\end{enumerate}
\eth
The special case of part (ii) when $\Theta$ is the set of all seeds was proved in
\cite[Theorem 4.6]{NTY2}. The proof there works for the general case.

The equality in (iv) is not proved in full generality in \cite[Lemma 3.3]{HLY}, but the proof is similar to the case of $R=\U_\ep(M_\ep, \B, \inv,\Theta)$ treated in \cite{HLY}. By the assumption on $R$, the root of unity mixed quantum torus-affine algebra 
$\TT_\ep(M_\ep)_\geq$ is a central localization of 
$R$:
\[
\TT_\ep(M'_\ep)_\geq = R[M'_\ep(\ell e_i)^{-1} ; i \in \ex].
\]
This implies at once the statement in part (iv). 

Since the algebras $\A_\ep(M_\ep,\B,\inv, \Theta) \subseteq \U_\ep(M_\ep,\B,\inv, \Theta)$
are subalgebras of $\TT_\ep(M_\ep)$, they are PI algebras \cite[Theorem A(1)]{NTY2} even without the coprimality assumption in \thref{cent-isom}. The last part (v) of the theorem makes this more precise in the upper quantum cluster algebra case. 

\bre{coprimality} \hfill
\begin{enumerate}
\item[(i)] Without loss of generality we can assume that the collection of seeds $\Theta$ of the root of unity quantum cluster algebra $\A_\ep(M_\ep, \B)$ contains the initial seed
$(M_\ep, \B)$. We will assume this throughout the paper for notational simplicity. 
\item[(ii)] In what follows, some results will assume the condition (Coprime) in \thref{cent-isom} and some will not. If a result does not mention this assumption, then it is stated in general.
\end{enumerate}
\ere
\subsection{Torus actions on $\A_\ep(M_\ep, \B, \inv, \Theta)$ and $\U_\ep(M_\ep, \B, \inv, \Theta)$}
\label{2.3}
Here we review torus actions on root of unity quantum cluster algebras following \cite[Sect. 3.3]{MNTY}. 

\bde{nB} For an exchange matrix $\B$, consider its nullity
\[
n(\B) := \dim \Ker(\B^\top)
\]
and the split torus (over $\KK$)
\begin{equation}
\label{TB}
T_{\B} := \Gset_m^{n(\B)}.
\end{equation}
\ede

For an exchange matrix $\B$ and $k\in \ex$, define an endomorphism $\mu_k$ of $\Zset^N$ by
\begin{equation}
\label{muKer}
\mu_k: \nu =(\nu_1, \dots, \nu_N) \in \Zset^N
\mt \nu':= (\nu_1, \dots,\nu_{k-1}, \nu_k', \nu_{k+1}, \dots, \nu_N), 
\end{equation}
where $\nu_k' := \nu\cdot[b^k]_+ - \nu_k$.
(Recall that $b^k$ denotes the $k$-th column of $\B$.) 
We have 
\begin{equation}
\label{muk-pm}
\nu \in \Ker(\B^\top)
\quad \Rightarrow \quad 
\nu_k' = \nu\cdot[b^k]_+ - \nu_k = 
\nu\cdot[b^k]_- - \nu_k \quad \mbox{and} \quad 
\nu' \in \Ker ( \mu_k (\B)^\top),
\end{equation}
see \cite[Lemma 2.3]{GSV} for the second fact. 

Each $\nu \in  \Ker(\B^\top)$ gives rise to a $\Gset_m$-action on the root of unity mixed quantum torus-affine algebra $\TT_\ep(M_\ep, \B)_\geq$ associated to the $\Zset$-grading  
\begin{equation}
\label{Zdeg}
\deg M_\ep(f):= \nu \cdot f, \; \; \forall f \in \Zset^N.
\end{equation}
Eq. \eqref{muk-pm} implies
\[
\deg \mu_k(M_\ep)(e_i) = \mu_k(\nu) \cdot e_i, \; \; \forall i \in [1,N].
\]
By recursive applications of this identity and using \eqref{muk-pm}, 
we obtain that $\A_\ep(M_\ep, \B, \inv, \Theta)$ and $\U_\ep(M_\ep, \B, \inv, \Theta)$
are graded subalgebras of $\TT_\ep(M_\ep, \B)_\geq$, i.e., they are preserved by the $\Gset_m$-action. 
This proves the first part of the following proposition:

\bpr{Caction} \cite[Sect. 3.3]{MNTY} Assume that $\ell$ is a positive integer satisfying condition (Coprime), $\ep^{1/2} \in \KK$ is a primitive $\ell$-th root of unity, and $\Theta$ is a connected collection of seeds of the root of unity quantum cluster algebra $\A_\ep(M_\ep, \B)$. 
\hfill
\begin{enumerate}
\item[(i)]
For every $\nu \in \Ker(\B^\top)$, there is a $\Gset_m$-action 
$\varphi_\nu$ on $\A_\ep(M_\ep, \B, \inv, \Theta)$ and $\U_\ep(M_\ep, \B, \inv, \Theta)$ 
associated to the $\Zset$-grading \eqref{Zdeg}. 
For every seed $(M'_\ep, \B') \in \Theta$ and $k \in [1,N]$, $M'_\ep(e_k)$ is homogeneous. 
\item[(ii)] Given a basis 
\[
\{\nu^1, \dots, \nu^{n(\B)} \}
\]
of the lattice $\Ker(\B^\top)$, 
we have an action of $T_{\B}$ on the algebras $\A_\ep(M_\ep, \B, \inv, \Theta)$ and
$\U_\ep(M_\ep, \B, \inv, \Theta)$ obtained by combining the $\Gset_m$-actions associated to $\nu^1$, $\dots,$ $\nu^{n(\B)}$. 
For every seed $(M'_\ep, \B') \in \Theta$ and $k \in [1,N]$, $M'_\ep(e_k)$ is homogeneous 
with respect to the corresponding $X(T_{\B})$-grading. 
\end{enumerate}
\epr

If $\ell$ satisfies condition (Coprime), then the torus action preserves the central subalgebras $\C\A_\ep(M_\ep, \B, \inv, \Theta)$ and $\C\U_\ep(M_\ep, \B, \inv, \Theta)$ and the quantum Frobenius map from \thref{cent-isom} (ii)-(iii) intertwines the two restrictions with the $\ell$-th power of the analogous torus actions on the algebras 
$\A(\B, \inv, \Theta)$ and $\U(\B, \inv, \Theta)$, see 
\cite[Lemma 3.4(ii)]{MNTY}.

\sectionnew{A discrete group of automorphisms of root of unity quantum cluster algebras}
\subsection{A faithful action of $\Ker \big( \, \ol{\B}^\top \big)$ on root of unity quantum cluster algebras} 
The endomorphism $\mu_k$ of $\Zset^N$ given by \eqref{muKer} descends to an endomorphism of $(\Zset/ \ell)^N$, which will be denoted by $\ol{\mu}_k$. It is given by 
\[
\ol{\mu}_k(\theta_1, \ldots, \theta_N) 
:=(\theta_1, \dots,\theta_{k-1}, \theta_k', \theta_{k+1}, \dots, \theta_N) \quad 
\mbox{with} \quad \theta_k' := \theta \cdot \ol{[b^k]_+} - \theta_k, 
\]
where $\theta := (\theta_1, \ldots, \theta_N) \in (\Zset/ \ell)^N$.

Each $\theta \in  \Ker \big( \, \ol{\B}^\top \big)$ gives rise to the following automorphism $(\theta \cdot)$ of the root of unity mixed quantum torus-affine algebra $\TT_\ep(M_\ep, \B)_\geq$
\begin{equation}
\label{Ker-act}
\theta \cdot M_\ep(f):= \ep^{\theta \cdot \ol{f}} M_\ep(f), \; \; \forall f \in \Zset^N.
\end{equation}
Let us extend it to an automorphism of the skew field of fractions $\FF$ of $\TT_\ep(M_\ep, \B)$. We have
\begin{align*}
\theta \cdot \mu_k(M_\ep)(e_k) & =\theta \cdot \Big( M_\ep([b^k]_+ - e_k) + M_\ep([b^k]_- - e_k)\Big)\\
&= \ep^{\theta \cdot (\ol{[b^k]_+} - \ol{e_k})} \Big( M_\ep([b^k]_+ - e_k) + \ep^{- \theta \cdot \ol{b^k}} M_\ep([b^k]_- - e_k) \Big)\\
&= \ep^{ \ol{\mu}_k(\theta) \cdot \ol{e_k}} \mu_k(M_\ep)(e_k),
\end{align*}
because $\theta \in  \Ker \big( \, \ol{\B}^\top \big)$ implies that 
$\theta \cdot \ol{b^k} =0$. Since $\mu_k(M_\ep)(e_i) = M_\ep(e_i)$ for $i \neq k$, we have 
$\theta \cdot \mu_k(M_\ep)(e_i) = \ep^{ \ol{\mu}_k(\theta) \cdot \ol{e_i}} \mu_k(M_\ep)(e_i)$, and hence, 
\[
\theta \cdot \mu_k(M_\ep)(f)= 
\ep^{ \ol{\mu}_k(\theta) \cdot \ol{f}} \mu_k(M_\ep)(f), 
\; \; \forall f \in \Zset^N.
\]
By recursive applications of this identity we obtain that this automorphism preserves the root of unity mixed quantum torus-affine algebra $\TT_\ep(M'_\ep, \B')_\geq$
associated to all seeds $(M'_\ep, \B') \sim (M_\ep, \B)$. 

Thus, for each $\theta \in \Ker \big( \, \ol{\B}^\top \big)$, eq. \eqref{Ker-act} defines an automorphism of $\A_\ep(M_\ep, \B, \inv, \Theta)$ and 
$\U_\ep(M_\ep, \B, \inv, \Theta)$. It is obvious that \eqref{Ker-act} defines an action of $\Ker \big( \, \ol{\B}^\top \big)$, i.e., $(\theta \cdot) (\theta' \cdot) = \big( (\theta + \theta') \cdot \big)$ for all $\theta, \theta' \in \Ker \big( \, \ol{\B}^\top \big)$. Furthermore, it is clear that the action fixes $\C\A_\ep(M_\ep, \B, \inv, \Theta)$ and $\C\U_\ep(M_\ep, \B, \inv, \Theta)$. If $\ell$ is odd (so that $\ep^n =1 \Rightarrow \ell \mid n$), the action is faithful because, for $\theta \in  \Ker \big( \, \ol{\B}^\top \big)$, 
\[
\theta \cdot M_\ep(e_i)= M_\ep(e_i), \; \forall i \in [1,N] \quad 
\Rightarrow \quad \theta \cdot \ol{e_i}=0, \; \forall i \in [1,N] \quad \Rightarrow \quad \theta = 0. 
\]
In summary, we obtain the following:

\bpr{ep-act} If $\ell$ is odd, 
for each connected collection $\Theta$ of seeds of the root of unity quantum cluster algebra $\A_\ep(M_\ep, \B)$, eq. \eqref{Ker-act} gives two injective group homomorphisms
\begin{align*}
&\Ker \big( \, \ol{\B}^\top \big) \hra \Aut \big(\A_\ep(M_\ep, \B, \inv, \Theta)/ \C\A_\ep(M_\ep, \B, \inv, \Theta) \big) \quad \mbox{and} \\
&\Ker \big( \, \ol{\B}^\top \big) \hra \Aut \big(\U_\ep(M_\ep, \B, \inv, \Theta)/ \C\U_\ep(M_\ep, \B, \inv, \Theta) \big).
\end{align*}
\epr
\bre{grad}
The above actions of the finite abelian group $\Ker \big( \, \ol{\B}^\top \big)$ on the partial root of unity quantum and upper quantum cluster algebras $\A_\ep(M_\ep, \B, \inv, \Theta)$ and $\U_\ep(M_\ep, \B, \inv, \Theta)$ are equivalent to gradings of these algebras by the character group of $\Ker \big( \, \ol{\B}^\top \big)$ (which is again a finite abelian group) with the property that every cluster/frozen variable is homogeneous. Such gradings on cluster algebras $\A(\B)$ were studied in \cite{WY} 
in many combinatorial situations related to tilings of the amplituhedron. 
\ere
\subsection{Reconstruction of automorphisms} 
Next, we obtain an effective characterization of the automorphisms in \prref{ep-act}, which will play a key role in the next section. First, recall the following:

\bde{nerve} \cite{F} A connected collection of seeds $\Theta$ of a classical or (root of unity) quantum cluster algebra is called a {\em{nerve}} if for each mutable direction $k \in \ex$, 
there are at least two seeds in $\Theta$ that are linked by the mutation $\mu_k$.
\ede

\bth{nerve} Assume that $\ell$ is an odd positive integer, $\ep^{1/2} \in \KK$ is a primitive $\ell$-th root of unity, $\Theta$ is a connected collection of seeds of the root of unity quantum cluster algebra $\A_\ep(M_\ep, \B)$ which is a nerve, and $\inv \subseteq \fr$ is arbitrary.  Let $R$ be an algebra that is sandwiched as 
\[
\A_\ep(M_\ep, \B, \inv, \Theta) \subseteq R \subseteq \U_\ep(M_\ep, \B, \inv, \Theta).
\]
If $\phi \in \Aut(R)$ is an automorphism such that for every seed $(M'_\ep, \B') \in \Theta$ and $i \in [1,N]$, 
\begin{equation}
\label{condII}
\phi( M'_\ep(e_i)) = \ep^m M'_\ep(e_i)
\end{equation}
for some $m \in \Zset / \ell$ depending on the seed and $i$, then $\phi$ is the restriction to $R$ of an automorphism \eqref{Ker-act} for some $\theta \in \Ker \big( \, \ol{\B}^\top \big)$.
\eth

\bre{notall} Given a subalgebra $R$ as in \thref{nerve}, not all automorphisms of $\U_\ep(M_\ep, \B, \inv, \Theta)$ of the form $(\theta \cdot)$ for $\theta \in \Ker \big( \, \ol{\B}^\top \big)$ will preserve $R$. The theorem claims that the automorphism $\phi$ is equal to an automorphism of the form $(\theta \cdot)$ that preserves $R$. 
\ere

For the proof of the theorem we will need the following lemma:
\ble{ind} Assume that $\B'$ and $\B''$ are two exchange matrices such that $\B'' = \mu_k(\B')$ for some $k \in \ex$. Denote by $(b')^i$ and $(b'')^i$ the $i$-th columns of $\B'$ and $\B''$ for $i \in \ex$. 
Let $\theta \in (\Zset/\ell)^N$ be such that 
$\ol{(b')^k}\cdot\theta =0$ and let $i \in \ex$ with $i \neq k$. Then 
$\ol{(b')^i} \cdot \theta=0$ if and only if $\ol{(b'')^i} \cdot \ol{\mu}_k(\theta)=0$.
\ele
\begin{proof} The lemma follows from the following identity:
\begin{align*}
\ol{(b'')^i} \cdot \ol{\mu}_k(\theta)
&= - \ol{b'_{ki}} (\theta \cdot \ol{[(b')^k]_+} - \theta_k)
+ \sum_{j\neq k}\left(\ol{b'_{ji}} + 
\frac{\ol{|b'_{jk}|b'_{ki}+b'_{jk}|b'_{ki}|}}{2}\right)\theta_j \\[6pt]
&= \ol{(b')^i} \cdot \theta + \ol{b'_{ki}}\left(\frac{\ol{[(b')^k]_+} + \ol{[(b')^k]_-}}{2} - \ol{[(b')^k]_+}\right)\cdot\theta + 
\frac{\ol{|b'_{ki}|}}{2}\, \ol{(b')^k}\cdot\theta \\[6pt]
&= \ol{(b')^i} \cdot \theta + \frac{\ol{|b'_{ki}|} - \ol{b'_{ki}}}{2}\,\ol{(b')^k}\cdot\theta 
= \ol{(b')^i} \cdot \theta.
\end{align*}
\end{proof}

\begin{proof}[Proof of \thref{nerve}] The assumption on $\phi$ gives that for every seed $\Sigma' = (M'_\ep, \B') \in \Theta$, there exists $\theta(\Sigma') \in (\Zset / \ell)^N$ such that
\[
\phi( M'_\ep(e_i)) = \ep^{\theta(\Sigma')_i} M'_\ep(e_i), \;\forall i \in [1,N].
\]
Denote $\theta:=\theta(\Sigma)$, where $\Sigma := (M_\ep, \B)$ is the initial seed, recall the convention in \reref{coprimality}(i). First we show that 

(*) If $\Sigma':=(M'_\ep, \B')$ and $\mu_k (\Sigma')$ are both elements of $\Theta$, then $\ol{(b')^k} \cdot \theta(\Sigma')=0$, where $(b')^k$ is the $k$-th column of $\B'$.

Indeed, for $\theta' := \theta(\Sigma')$ we have 
\begin{align*}
&\ep^{\big( \theta(\mu_k(\Sigma')) \big)_k} 
\big( M'_\ep([(b')^k]_+ - e_k) + M'_\ep([(b')^k]_- - e_k) \big)
= \phi (\mu_k (M'_\ep)(e_k)) 
\\
= 
&\ep^{\theta' \cdot (\ol{[(b')^k]_+} - \ol{e_k})}
M'_\ep([(b')^k]_+ - e_k) + 
\ep^{\theta' \cdot (\ol{[(b')^k]_-} - \ol{e_k})}
M'_\ep([(b')^k]_- - e_k).
\end{align*}
Therefore, $\theta' \cdot (\ol{[(b')^k]_+} - \ol{e_k}) = \theta' \cdot (\ol{[(b')^k]_-} - \ol{e_k})$, which proves (*). 

\leref{ind}, the property (*) and the assumption that the collection of seeds $\Theta$ is a nerve imply that 
\[
\theta \cdot \ol{b^k} =0, \; \forall \; k \in \ex.
\]
Hence, $\theta \in \Ker \big( \, \ol{\B}^\top \big)$. This implies that 
\[
\phi|_{\A_\ep(M_\ep, \B, \inv, \Theta)} = 
(\theta \cdot)|_{\A_\ep(M_\ep, \B, \inv, \Theta)},
\]
and therefore the extensions of $\phi$ and $(\theta \cdot)$ to the skew field of fractions $\FF$ of $\A_\ep(M_\ep, \B)$ coincide. Restricting them to $R$ gives that $\phi$ equals the restriction of $(\theta \cdot)$ to $R$. 
\end{proof}
\sectionnew{Relative and ozone automorphism groups for sandwiched subalgebras}
\lb{results}
\subsection{Relative automorphisms for sandwiched algebras of quantum tori}
Throughout this subsection, for a positive integer $n$, denote by $n'$ the (unique) positive integer defined as follows. If $\charr \KK =0$, set 
\begin{equation}
    \label{nn'1}
n' := n.
\end{equation}
If $\charr \KK = p$, let
\begin{equation}
    \label{nn'2}
n = n' p^k, \quad \mbox{where} \quad p \nmid n', \; k \in \Nset.
\end{equation}
Denote by 
\begin{equation}
\label{mun}
C_{n'} := \{ c \in \KK \mid c^{n'} = 1 \} 
\end{equation}
the group of $n'$-th roots of unity in $\KK$.
(Often this group is denoted by $\mu_{n'}$ but this notation interferes with the notation for mutation in cluster algebras.) By the definition of $n'$, 
\[
C_{n'} = \{ c \in \KK \mid c^n = 1 \}.
\]
If $\KK$ contains a primitive $n'$-th root of unity $\zeta$ (this can always be achieved in a finite extension of $\KK$), then  
\[
C_{n'} = \{ \zeta^k \mid k \in \Zset/n' \} \cong \Zset/ n'.
\]

Given a matrix ${\mathbf{q}}:=(q_{ij})_{i,j=1}^N \in M_N(\KK^\times)$, which is multiplicatively skew-symmetric, that is,  
$q_{ii} =1$ and 
$q_{ij} = q_{ji}^{-1}$ for all $i, j \in [1,N]$, we consider the quantum torus 
\[
T_{\mathbf{q}}:= \frac{\KK \langle x_1^{\pm 1}, \ldots,  x_N^{\pm 1} \rangle}{(x_i x_j -q_{ij} x_j x_i  :  1 \leq j < i \leq N)}.
\]
Denote by $A_{\mathbf{q}}$ the corresponding quantum affine space algebra, that is, the unital subalgebra generated by $x_1, \ldots, x_N$. For $\nu:=(\nu_1, \ldots, \nu_N) \in \Zset^N$ denote 
\[
X^{[\nu]}:=x_1^{\nu_1} \cdots x_N^{\nu_N}.
\]
We have 
\begin{equation}
\label{mult-non-norm}
X^{[\pi]} X^{[\nu]} = c(\pi, \nu) X^{[\pi + \nu]}, \quad \forall \pi, \nu \in \Zset^N
\end{equation}
for some $c(\pi,\nu) \in \KK^\times$.
\medskip

For a positive integer $n$, denote by $A_{\mathbf{q}}^n$ and 
$T_{\mathbf{q}}^n$ the subalgebras of 
$A_{\mathbf{q}}$ and $T_{\mathbf{q}}$ generated by 
$x_1^n, \ldots, x_N^n$ and $x_1^{\pm n}, \ldots, x_N^{\pm n}$, respectively. 

\bth{gen} Consider a quantum torus $T_{\mathbf{q}}$ and a positive integer $n$. Let $R$ be an algebra 
sandwiched as
\[
A_{\mathbf{q}} \subseteq R \subseteq T_{\mathbf{q}}.
\]
Then, 
\[
\End(R/A_{\mathbf{q}}^n) = \Aut(R/A_{\mathbf{q}}^n)
\]
and for all $\phi \in \Aut(R/A_{\mathbf{q}}^n)$ there exist 
$\zeta_1, \ldots, \zeta_N \in C_{n'}$ such that 
\[
\phi(x_i) = \zeta_i x_i, \quad \forall i \in [1,N].
\]
\eth

\bco{repeat} For an algebra $R$ such that $A_{\mathbf{q}} \subseteq R \subseteq T_{\mathbf{q}}$ and a subalgebra $S \subseteq R$ such that 
 $A_{\mathbf{q}}^n \subseteq S \subseteq T_{\mathbf{q}}^n$, we have:
 \begin{enumerate}
 \item[(i)] $\End(R/S)=\Aut(R/S)$ and 
 \item[(ii)]
 $\Aut(R/S)$ is isomorphic to a subgroup of $(\Zset/n')^N$ and consists of the automorphisms of $T_{\mathbf{q}}$ of the form 
\begin{equation}
\label{zeta}
\phi(x_i) = \zeta^{\theta_i} x_i, \quad \forall i \in [1,N]
\end{equation} 
that preserve $R$, where $\zeta \in \KK$ is a (fixed) primitive $n'$-th root of unity and $(\theta_1, \ldots, \theta_N) \in (\Zset/n')^N$.
\end{enumerate}
\eco
\begin{proof}(i) Since $A_{\mathbf{q}}^n \subseteq S \subseteq R$, we have
$\End(R/S) \subseteq \End(R/A_{\mathbf{q}}^n) =\Aut(R/A_{\mathbf{q}}^n)$. So, by \thref{gen}, every $\phi \in \End(R/S)$ is an automorphism of $R$, and it fixes $S$ by assumption; that is, $\End(R/S) = \Aut(R/S)$.

(ii) By \thref{gen} and $\Aut(R/S) \subseteq \Aut(R/A_{\mathbf{q}}^n)$, every $\phi \in \Aut(R/S)$ has the form \eqref{zeta}, and $\phi \mapsto (\theta_1, \ldots, \theta_N)$ is an injective homomorphism into $(\Zset/n')^N$; conversely, an automorphism of the form \eqref{zeta} that preserves $R$ fixes $T_{\mathbf{q}}^n \supseteq S$, so it lies in $\End(R/S) = \Aut(R/S)$ by part (i).
\end{proof}
We will need the following auxiliary result for the proof of \thref{gen}. 

\bpr{main} For each quantum torus $T_{\mathbf{q}}$ and $\nu \in \Zset^N$, the 
solutions of the equation 
\[
a^n =(X^{[\nu]})^n
\]
in $T_{\mathbf{q}}$ are 
\[
a=cX^{[\nu]} \quad \mbox{for} \quad c \in C_{n'}.
\]
\epr
\begin{proof}
Write
\[
a=\sum_{h\in \Delta} c(h) X^{[h]},
\]
where $\Delta \subset \mathbb Z^N$ is a finite set and $c(h) \in \KK^\times$. Let $\Psi$ be the convex hull of $\Delta$. We will show that every vertex of $\Psi$ contributes to a monomial of $a^n$ in a unique fashion, from which we will deduce that $\Psi$ has only one vertex. Let $h_0$ be a vertex of $\Psi$. We can choose a linear map $\lambda : \mathbb R^N \rightarrow \mathbb R$ such that $\lambda(h_0)>\lambda (h)$ for all $h \in \Psi \backslash \{h_0 \}$, and in particular for all
$h\in \Delta \backslash \{h_0\}$. Now consider the terms in the expansion of $a^n$. Each term is of the form 
\[
c(h_1) \cdots c(h_n) X^{[h_1]} \cdots X^{[h_n]}
\]
for some $h_1, \ldots, h_n \in \Delta$ and 
is a scalar multiple of $X^{[h_1 + \cdots + h_n]}$ because of \eqref{mult-non-norm}.
In order for such a term to be a scalar multiple of $X^{[n h_0]}$,
we must have 
$$h_1+\dots + h_n = n h_0$$ 
for some $h_i \in \Delta$, $1\leq i \leq n$. 
Therefore,
$$\lambda(h_1)+\dots +\lambda(h_n) = n \lambda(h_0).$$
Since $h_i \in \Psi$, we have $\lambda(h_i) < \lambda (h_0)$ for $h_i \neq h_0$. Thus, 
$$h_1=\cdots=h_n = h_0.$$ 
Therefore, $\Psi$ must have only one vertex, and thus, $\Delta = \{ h_0 \}$. So, $a = c(h_0) X^{[h_0]}$ and 
$$c(h_0)^n(X^{[h_0]})^n =(X^{[\nu]})^n,$$
which is equivalent to saying that $h_0 = \nu$ and $c(h_0)^n =1$. The last equation is equivalent to $c(h_0) \in C_{n'}$.
\end{proof}
\begin{proof}[Proof of \thref{gen}] Let $\phi \in \End(R/A_{\mathbf{q}}^n)$. 
Since $x_i^n \in A_{\mathbf{q}}^n$, 
$\phi(x_i)^n=\phi(x_i^n) = x_i^n$. \prref{main} implies that 
$$\phi(x_i)=c_i x_i$$
for some $c_i \in C_{n'}$. Such a map $\phi$ extends to an automorphism of  $T_{\mathbf{q}}$ of the form \eqref{zeta} satisfying $\phi(R) \subseteq R$.
In particular, $\Ker \phi =\{0\}$. Since $\phi^{n'} = \id$, 
\[
T_{\mathbf{q}} = \bigoplus_{m \in \Zset/n'} (T_{\mathbf{q}})_{\zeta^m},
\]
where $(T_{\mathbf{q}})_{\zeta^m}$ 
are the eigenspaces of $\phi$ acting on $T_{\mathbf{q}}$ with eigenvalues $\zeta^m$ for $m \in \Zset/n'$.
Because $\phi(R) \subseteq R$, we have 
\[
R = \bigoplus_{m \in \Zset/n'} R \cap (T_{\mathbf{q}})_{\zeta^m},
\]
which implies that $\Im \phi  =R$. Therefore,  $\phi \in \Aut(R/A_{\mathbf{q}}^n)$.
\end{proof}
\subsection{Relative automorphisms for sandwiched algebras of root of unity quantum cluster algebras}
\bth{aut-qcl} Assume that $\ell$ is a positive odd integer coprime to $\charr \KK$, $\ep^{1/2} \in \KK$ is a primitive $\ell$-th root of unity, $\Theta$ is a connected collection of seeds of the root of unity quantum cluster algebra $\A_\ep(M_\ep, \B)$ that is a nerve, and $\inv \subseteq \fr$ is arbitrary. 
Consider a sandwiched subalgebra  
\[
\A_\ep(M_\ep, \B, \inv, \Theta) \subseteq R \subseteq \U_\ep(M_\ep, \B, \inv, \Theta).
\]
\begin{enumerate}
\item[(i)] For every subalgebra $C$ of $R$ sandwiched as 
\[
\C\A_\ep(M_\ep, \B, \inv, \Theta) \subseteq C \subseteq \C\U_\ep(M_\ep, \B, \inv, \Theta), 
\]
(which is necessarily central in $R$), we have
\begin{equation}
\label{EndAut}
\End(R/C) = \Aut(R/C).
\end{equation}
Furthermore, the relative automorphism group $\Aut(R/C)$ is isomorphic to a subgroup of $\Ker \big( \, \ol{\B}^\top \big)$ and consists of the automorphisms of $\U_\ep(M_\ep, \B, \inv, \Theta)$ of the form $(\theta \cdot)$ that preserve $R$ (recall \eqref{Ker-act}).
\item[(ii)] The ozone group $\Oz(R)$ of $R$ is isomorphic to a subgroup of 
\[
\Ker \big( \, \ol{\B}^\top \big) \cap \Big( 
\ol{\Ker \Om_{M_\ep}} \Big)^\perp
\]
and consists of the automorphisms of $\U_\ep(M_\ep, \B, \inv, \Theta)$ of the form $(\theta \cdot)$ that preserve $R$. The orthogonal complement is taken with respect to the dot product \eqref{dot-prod} on $(\Zset/\ell)^N$.
\end{enumerate}
\eth
Recall from Sect. \ref{2.2} that, for a root of unity toric frame $M_\ep$, $\Om_{M_\ep}$ denotes both the additive bicharacter $\Om_{M_\ep} : \Zset^N \times \Zset^N \to \Zset/\ell$ and its matrix in $M_N(\Zset/\ell)$. In part (ii) of the theorem we use the former convention.

\begin{proof} (i) The equality \eqref{EndAut} follows from \coref{repeat}(i). 
Let $\phi \in \Aut(R/C)$. For every seed $(M'_\ep, \B') \in \Theta$ and $i \in [1,N]$, we have
\[
\phi(M'_\ep(e_i))^\ell = \phi(M'_\ep( e_i)^\ell) = M'_\ep(e_i)^\ell,
\] 
so \prref{main} gives $\phi(M'_\ep(e_i)) = \zeta_i\, M'_\ep(e_i)$ for some $\zeta_1, \ldots, \zeta_N \in C_\ell$. \thref{nerve} implies that $\phi = (\theta \cdot)$ for some 
$\theta \in \Ker \big( \, \ol{\B}^\top \big)$, and as a consequence, $\theta \cdot R = R$. 

In the opposite direction, assume that $\theta \in \Ker \big( \, \ol{\B}^\top \big)$ is such that $\theta \cdot R \subseteq R$. By \prref{ep-act}, $(\theta \cdot)$ fixes
$\C\U_\ep(M_\ep, \B, \inv, \Theta) \supseteq C$, so $(\theta \cdot) \in \End(R/C)$,
and \eqref{EndAut} gives $(\theta \cdot) \in \Aut(R/C)$.

(ii) Assume that $\phi \in \Oz(R)$. Part (i) implies that 
$\phi = (\theta \cdot)$ for some 
$\theta \in \Ker \big( \, \ol{\B}^\top \big)$ such that 
$\theta \cdot R = R$. 
Take $f \in \Ker \Om_{M_\ep} \subset \Zset^N$. Choose $n_i \in \Nset$, $i \in [1,N]$, such that
\[
g := f + \ell \sum_{i=1}^N n_i e_i \in 
\Ker \Om_{M_\ep} \cap \Nset^N.
\]
By \eqref{Tep}, $y := M_\ep(g)$ is a nonzero scalar
multiple of a product of the elements $M_\ep(e_i)$, so
$y \in \A_\ep(M_\ep, \B, \inv, \Theta) \subseteq R$, and by \eqref{ZT},
$y \in Z(\TT_\ep(M_\ep)_\geq)$. Since $R \subseteq \TT_\ep(M_\ep)_\geq$, it follows that
$y \in Z(R)$. Thus,
\[
y = \phi(y) = \ep^{\theta \cdot \ol{f}} y.
\]
Therefore, $\theta \in \Big( \ol{\Ker \Om_{M_\ep}} \Big)^\perp$. This proves that 
\[
\phi = (\theta \cdot) \quad 
\mbox{for} \quad 
\theta \in \Ker \big( \, \ol{\B}^\top \big) \cap \Big( 
\ol{\Ker \Om_{M_\ep}} \Big)^\perp \quad \mbox{such that} \quad \theta \cdot R \subseteq R.
\]
In the opposite direction, if $\phi$ is an automorphism 
of $\TT_\ep(M_\ep)_\geq$ of this form, then 
$\phi$ restricts to an endomorphism of $R$, which has to be an automorphism by \coref{repeat}(i). Furthermore, if $z \in Z(R)$, then $z \in Z(\TT_\ep(M_\ep)_\geq)$ and by \eqref{ZT}, $z$ is given by a finite sum of the form 
\[
z = \sum_{f \in \Ker \Om_{M_\ep}} c(f) M_\ep(f). 
\]
Hence, 
\[
\phi(z) = \sum_{f \in \Ker \Om_{M_\ep}} 
\ep^{\theta \cdot \ol{f}}
c(f) M_\ep(f) = z.
\]
So, $\phi \in \Oz(R)$.
\end{proof}

As an immediate consequence of \thref{aut-qcl}, we have:

\bco{Oz} Assume that $\KK$, $\ell$, $\ep^{1/2}$, $\B$, $\Theta$ and $\inv$ are as in \thref{aut-qcl}.
\begin{enumerate}
\item[(i)] The Frobenius endomorphism monoids and automorphism groups of the corresponding partial root of unity  quantum and upper quantum cluster algebras are given by 
\begin{align*}
&\End \big( \A_\ep(M_\ep, \B, \inv, \Theta)/
\C\A_\ep(M_\ep, \B, \inv, \Theta) \big) = \\
& \Aut \big( \A_\ep(M_\ep, \B, \inv, \Theta)/
\C\A_\ep(M_\ep, \B, \inv, \Theta) \big) = \Ker \big( \, \ol{\B}^\top \big),
\end{align*}
and 
\begin{align*}
&\End \big( \U_\ep(M_\ep, \B, \inv, \Theta)/
\C\U_\ep(M_\ep, \B, \inv, \Theta) \big) = \\
& \Aut \big( \U_\ep(M_\ep, \B, \inv, \Theta)/
\C\U_\ep(M_\ep, \B, \inv, \Theta) \big) = \Ker \big( \, \ol{\B}^\top \big),
\end{align*}
where the equalities mean that every automorphism is realized
as the action $(\theta \cdot)$ of an element of the corresponding group in the right-hand side.
\item[(ii)] The ozone groups of the partial root of unity quantum and upper quantum cluster algebras are given by 
\[
\Oz \big( \A_\ep(M_\ep, \B, \inv, \Theta) \big) =
\Oz \big( \U_\ep(M_\ep, \B, \inv, \Theta) \big) = 
\Ker \big( \, \ol{\B}^\top \big) \cap \big( \ol{\Ker \Om_{M_\ep}} \big)^\perp.
\]
\end{enumerate}
\eco
\bex{G2-noncoprime} 
Let $\KK$ be a field with $\charr \KK \neq 3$, $\ell = 9$, and let
$\ep^{1/2} \in \KK$ be a primitive $9$th root of unity. Let $\ex = \{1,2\}$, $\fr = \{3\}$ and
\[
\widetilde{B} = \begin{bmatrix} 0 & 1 \\ -3 & 0 \\ 3 & 0 \end{bmatrix}, \quad
\Omega = \begin{bmatrix} 0 & 1 & 0 \\ -1 & 0 & 0 \\ 0 & 0 & 0 \end{bmatrix}, \quad
D = \begin{bmatrix} 3 & 0 \\ 0 & 1 \end{bmatrix}.
\]
The principal part of $\widetilde{B}$ is of type $G_2$, it is skew-symmetrized by $D$, and
$\widetilde{B}^\top \Omega = [D \; 0]$. The assumptions of \thref{aut-qcl} are satisfied, but condition (Coprime) fails since $3 \mid \ell$. For $\theta = (\theta_1,\theta_2,\theta_3)^{\top}\in(\mathbb Z/9)^3$, we have $\ol{\widetilde{B}}^\top \theta = (3(\theta_3 - \theta_2), \theta_1)$, so \thref{aut-qcl} gives
\begin{align*}
    \Aut \big( \R_\ep(M_\ep, \widetilde{B}, \inv, \Theta)/\C\R_\ep(M_\ep, \widetilde{B}, \inv, \Theta) \big) &
\cong \{ (0, \theta_2, \theta_3) \mid \theta_3 \equiv \theta_2 \mod 3 \}\\
& \cong \Zset/9 \oplus \Zset/3
\end{align*}
which is not isomorphic to $(\Zset/9)^{|\fr|}$. Viewing $\Omega$ as a bicharacter $\Zset^3 \times \Zset^3 \to \Zset/9$, we have
$\Ker \Om_{M_\ep} = 9\Zset \oplus 9\Zset \oplus \Zset$, so $\ol{\Ker \Om_{M_\ep}} = \{ (0,0,c) \mid c \in \Zset/9 \}$. \coref{Oz}(ii) gives
\[
\Oz \big( \A_\ep(M_\ep, \widetilde{B}, \inv, \Theta) \big)
= \Oz \big( \U_\ep(M_\ep, \widetilde{B}, \inv, \Theta) \big)
\cong \{ (0, \theta_2, 0) \mid \theta_2 \in \{0, 3, 6\} \} \cong \Zset/3.
\]
The assumptions of \coref{Oz} are weaker than condition (Coprime). The above example shows that in this generality the Frobenius automorphism groups need not be isomorphic to $(\Zset/\ell)^{|\fr|}$ and the conclusions of \leref{nullity}(ii), \thref{torus}, \coref{isom-ker} and
\thref{normal-Oz} in the following section fail without condition (Coprime).
\eex

\sectionnew{Frobenius and ozone automorphism groups via normal elements and torus actions}

In this section we prove two other explicit theorems for the Frobenius and ozone automorphism groups of cluster algebras. The first one concerns both groups and is based on actions of split tori on root of unity quantum cluster algebras. The second one concerns only ozone groups and uses normal elements. 

\subsection{Frobenius and ozone automorphism groups of cluster algebras via torus actions}
\ble{nullity} If $\ell$ is a positive integer satisfying condition (Coprime) for an exchange matrix $\B$, 
then the following hold:
\begin{enumerate}
\item[(i)] The exchange matrix $\B$ has full rank, i.e., its nullity $n(\B)$ equals $|\fr|$, recall \deref{nB}.  
\item[(ii)] The order of the abelian group $\Ker \big( \, \ol{\B}^\top \big)$ equals $\ell^{|\fr|}$.
\end{enumerate}
\ele
\begin{proof} 
(i) By \eqref{compatible}, we have 
\[
\ol{\B}^\top \Omega  = 
\begin{bmatrix}
\, \ol{D} \; 0
\end{bmatrix} . 
\]  
The diagonal entries of the diagonal matrix $\ol{D}$ equal
$\ol{d_i}$ for $i \in \ex$ and are units of $\Zset/ \ell$ because $\ell$ and $d_i$ are relatively prime by condition (Coprime). So, $\ol{\B}$ has (full) rank $|\ex|$ and thus, $\B$ also has (full) rank $|\ex|$. Therefore, its nullity is 
\[
n(\B) = N - |\ex| = |\fr|.
\]

(ii) Consider $\ol{\B}^\top$ as a homomorphism from 
$(\Zset/\ell)^N$ to $(\Zset/\ell)^{|\ex|}$. It is surjective since \eqref{compatible} implies that $\ol{\B}^\top ( \ol{d_i}^{-1} \Omega \ol{e_i}) = \ol{e_i}$ for all $i \in \ex$. Therefore, 
\[
\Big| \Ker \big( \, \ol{\B}^\top \big) \Big| = \frac{|\Zset/\ell|^N}{\Big| \Im \big( \, \ol{\B}^\top \big) \Big|} = 
\frac{\ell^N}{\ell^{|\ex|}}= \ell^{|\fr|}.
\]
\end{proof}

For a split torus $T = \Gset_m^n$ over $\KK$ of rank $n$, its group of $\KK$-points is $T(\KK)\cong (\KK^\times)^n$. Its $\ell$-torsion subgroup will be denoted by 
\[
T(\KK)[\ell] = \{ y \in (\KK^\times)^n \mid y^\ell = 1 \} = C_\ell^n, 
\]
where in the last equality we use the assumption that $\KK$ contains a primitive $\ell$-th root of unity, namely $\ep^{1/2}$.

\bth{torus} Assume that $\ell$ is a positive integer satisfying condition (Coprime), $\ep^{1/2} \in \KK$ is a primitive $\ell$-th root of unity, 
$\Theta$ is a connected collection of seeds of the root of unity quantum cluster algebra $\A_\ep(M_\ep, \B)$ that is a nerve, and $\inv \subseteq \fr$ is arbitrary. 
\begin{enumerate}
\item[(i)] Assume that $T$ is a split torus of rank 
$|\fr|$ acting on $\A_\ep(M_\ep, \B, \inv, \Theta)$ (respectively $\U_\ep(M_\ep, \B, \inv, \Theta)$)
such that the $\ell$-torsion subgroup $T(\KK)[\ell]$ of the group 
of $\KK$-points $T(\KK) \cong (\KK^\times)^{|\fr|}$ acts faithfully and the cluster variables $M'_\ep(e_i)$ are 
$T(\KK)$-eigenvectors for all $(M'_\ep, \B') \in \Theta$ and $i \in [1,N]$.

The Frobenius automorphism groups of these algebras are given by 
\begin{align*}
&\Aut \big(
\A_\ep(M_\ep, \B, \inv, \Theta)/
\C\A_\ep(M_\ep, \B, \inv, \Theta) \big) = \\
&\Aut \big(
\U_\ep(M_\ep, \B, \inv, \Theta)/
\C\U_\ep(M_\ep, \B, \inv, \Theta) \big) = 
T(\KK)[\ell] = C_\ell^{|\fr|} \cong (\Zset/\ell)^{|\fr|},
\end{align*}
recall \eqref{mun}. Here the equalities mean that every element of the Frobenius automorphism 
groups comes from the action of the $\ell$-torsion subgroup of $T$.

The ozone groups of these algebras are given by 
\begin{align*}
&\Oz \big( \A_\ep(M_\ep, \B, \inv, \Theta) \big) =
\Oz \big( \U_\ep(M_\ep, \B, \inv, \Theta) \big) = 
\\
&\{ t \in C_\ell^{|\fr|} \mid t \cdot M_\ep(f) = M_\ep(f), \; \forall 
f \in \Ker (\Om_{M_\ep})\}.
\end{align*}
\item[(ii)] There is a canonical torus with the properties in the first part, namely the torus $T_{\B}$ from \eqref{TB} with 
the action in \prref{Caction} constructed from a basis of the lattice $\Ker(\B^\top) \subseteq \Zset^N$.  
\end{enumerate}
\eth

For all $t \in C_\ell^{|\fr|}$ and $t \in T_{\B}(\KK)[\ell]$, the automorphism $(t \cdot)$ of $\A_\ep(M_\ep, \B, \inv, \Theta)$ and $\U_\ep(M_\ep, \B, \inv, \Theta)$ extends to an automorphism of $\TT_\ep(M_\ep)$ (which is a localization of $\A_\ep(M_\ep, \B, \inv, \Theta)$ 
and $\U_\ep(M_\ep, \B, \inv, \Theta)$ by the multiplicative subset generated by the elements in the Frobenius center  
$M_\ep(\ell e_i) \in \C\A_\ep(M_\ep, \B, \inv, \Theta)$
for $i \in [1,N]$) by setting $t \cdot M_\ep(\ell e_i)^{-1} := M_\ep(\ell e_i)^{-1}$.
These are the actions on $M_\ep(f)$ that are used in both parts of the theorem. 


Note that \leref{nullity}(i) gives that the rank of the split torus $T_{\B}$ equals $|\fr|$.

\begin{proof} \coref{Oz} and \leref{nullity}(ii) imply that the orders of both groups 
\begin{align*}
&\Aut\big(\A_\ep(M_\ep, \B, \inv, \Theta)/
\C\A_\ep(M_\ep, \B, \inv, \Theta) \big) \quad \mbox{and}
\\
&\Aut \big(
\U_\ep(M_\ep, \B, \inv, \Theta)/
\C\U_\ep(M_\ep, \B, \inv, \Theta) \big)
\end{align*}
equal $\ell^{|\fr|}$.
Next we show that there are injective homomorphisms from the abelian group $T(\KK)[\ell] = C_\ell^{|\fr|}$ of order $\ell^{|\fr|}$ to both groups, which proves that the groups coincide with the images of the two homomorphisms. This will prove the statement about relative automorphism groups in part (i). If $t \in T(\KK)[\ell]$, then by the assumption about the $T$-action,
\[
t \cdot M'_\ep(e_i) = \zeta M'_\ep(e_i), \quad 
\forall (M'_\ep, \B') \in \Theta, \; i \in [1,N]
\]
for some $\zeta \in \KK^\times$ depending on the seed 
$(M'_\ep, \B')$ and $i \in [1,N]$.
Furthermore, since $t \in T(\KK)[\ell]$, 
\[
M'_\ep(e_i) = t^\ell \cdot M'_\ep(e_i) = \zeta^\ell M'_\ep(e_i).
\]
Therefore $\zeta^\ell =1$ and
\[
t \cdot M'_\ep(e_i)^\ell  = \zeta^\ell M'_\ep(e_i)^\ell = M'_\ep(e_i)^\ell,
\]
which implies that $(t \cdot)$ fixes 
$\C\A_\ep(M_\ep, \B, \inv, \Theta)$
(respectively $\C\U_\ep(M_\ep, \B, \inv, \Theta)$). This completes the proof of the statement about the Frobenius automorphism groups in part (i). 

Next, we show the statement about the ozone group of  
$\A_\ep(M_\ep, \B, \inv, \Theta)$; the statement about 
$\Oz(\U_\ep(M_\ep, \B, \inv, \Theta))$ is proved in a similar manner. Let 
\begin{align}
\nn
\phi \in &\Oz(\A_\ep(M_\ep, \B, \inv, \Theta)) \subseteq 
\\
&\Aut\big(\A_\ep(M_\ep, \B, \inv, \Theta)/\C\A_\ep(M_\ep, \B, \inv, \Theta) \big).
\label{in}
\end{align}
By the statement about relative automorphism groups in part (i), $\phi = (t \cdot)$ for some $t \in T(\KK)[\ell]$. 
As noted after the statement of the theorem, $(t \cdot)$ extends to an automorphism of 
\[
\TT_\ep(M_\ep) = \A_\ep(M_\ep, \B, \inv, \Theta)
[M_\ep(\ell e_i)^{-1}, i \in [1,N]]
\]
by setting 
$t \cdot M_\ep(\ell e_i)^{-1} = c^{-1} M_\ep(\ell e_i)^{-1}$
for the scalar $c \in \KK^\times$ such that 
$t \cdot M_\ep(\ell e_i) = c M_\ep(\ell e_i)$. 
Let $f \in \Ker \Om_{M_\ep}$ with $(f)_i \geq 0$ for
$i \in \ninv$.  There exist
$n_i \in \Nset$ for $i \in \ex \sqcup \inv$ such that
\[
f + \ell \sum_{i \in \ex \sqcup \inv} n_i e_i \in 
\Ker \Om_{M_\ep} \cap \Nset^N.
\]
The element
\[
y:=M_\ep(f) \prod_{i \in \ex \sqcup \inv} M_\ep(n_i e_i)^{\ell} 
\]
belongs to $Z(\A_\ep(M_\ep, \B, \inv, \Theta))$ because it is a scalar multiple of $M_\ep(f + \sum_{i \in \ex \sqcup \inv} \ell n_i e_i)$ by \eqref{Tep} and because of \eqref{ZT} and \thref{cent-isom}(iv). 
This implies $t \cdot y = y$. 
Furthermore, $t \cdot M_\ep(n_i e_i)^{\ell} = M_\ep(n_i e_i)^{\ell}$ for all $i \in [1,N]$ because $t \in T(\KK)[\ell]$. Hence $t \cdot M_\ep(f) = M_\ep(f)$. 

It remains to show that, for every 
$t \in C_\ell^{|\fr|}$ such that $t \cdot M_\ep(f) = M_\ep(f)$ for all $f \in \Ker (\Om_{M_\ep})$ with $(f)_i \geq 0$ for $i \in \ninv$, we have $(t \cdot) \in \Oz(\A_\ep(M_\ep, \B, \inv, \Theta))$. For this, we need to show that for all $z \in Z(\A_\ep(M_\ep, \B, \inv, \Theta))$, $t \cdot z = z$. By \eqref{ZT} and \thref{cent-isom}(iv), $z$ is a linear combination of terms of the form $M_\ep(f)$ for $f \in \Ker (\Om_{M_\ep})$ with $(f)_i \geq 0$ for $i \in \ninv$. Since all of them are fixed by $(t \cdot)$, the same is true for $z$ as well.

(ii) Under the condition (Coprime), the rank of the torus $T_{\B}$, defined in \eqref{TB}, equals $|\fr|$ by \leref{nullity}(i). We will show that the group of $\KK$-points $T_{\B}(\KK) \cong (\KK^\times)^{|\fr|}$ of the torus $T_{\B}$ acts faithfully on $\A_\ep(M_\ep, \B, \inv, \Theta)$. This will imply that the statement in part (ii) follows from the one in part (i); for this implication we need the weaker statement that only the $\ell$-torsion subgroup $T_{\B}(\KK)[\ell]$ acts 
faithfully on $\A_\ep(M_\ep, \B, \inv, \Theta)$.

By \leref{nullity}(i), $\B$ has full rank (equal to $|\ex|$). Choose matrices $A_1 \in \mathrm{GL}_{\ex}(\Zset)$ and 
$A_2 \in \mathrm{GL}_N(\Zset)$ such that 
$A_1 \B^\top A_2$ is the Smith normal form of $\B^\top$ with zero columns in positions $1, \ldots, |\fr|$. Then 
\[
\Ker \B^\top = \Zset A_2 e_1 \oplus \cdots \oplus \Zset A_2 e_{|\fr|}.
\]
Choose $\{\nu^i := A_2 e_i \mid i \in [1, |\fr|]\}$ as a basis of the lattice $\Ker \B^\top$ as in \prref{Caction}(ii). Let $j \in [1, |\fr|]$.
For the corresponding action of $T_{\B} \cong \Gset_m^{N - |\ex|}$ on the root of unity quantum torus $\TT_\ep(M_\ep)$ we have 
\begin{align*}
\deg M_\ep((A_2^{-1}){^\top} e_j ) &= 
(A_2 e_1 \cdot (A_2^{-1}){^\top} e_j, \dots, A_2e_{N-|\ex|} \cdot (A_2^{-1})^{\top} e_j) 
\\
&=(e_1 \cdot e_j, \ldots, e_{N- |\ex|} \cdot e_j) =
e_j.
\end{align*}
This means that the group generated by the characters of $T_\B$ for the nonzero components for its action on $\TT_\ep(M_\ep)$ 
is the full character lattice $X(T_\B)$. 
Therefore the action of $T_\B(\KK)$ on 
$\TT_\ep(M_\ep)$ is faithful. Since the latter is a localization of $\A_\ep(M_\ep, \B, \inv, \Theta)$, the action of $T_\B$ on $\A_\ep(M_\ep, \B, \inv, \Theta)$ is also faithful, which completes the proof of part (ii). 
\end{proof}
Combining \coref{Oz} and \thref{torus}(ii), we obtain the following:
\bco{isom-ker} If $\ell$ is a positive integer satisfying condition (Coprime) for an exchange matrix $\B$, then 
\[
\Ker \big( \, \ol{\B}^\top \big) \cong (\Zset/ \ell)^{|\fr|}.
\]
\eco

\subsection{Ozone groups of cluster algebras via normal elements}

Next we prove a classification theorem for the ozone groups of root of unity quantum and upper quantum cluster algebras 
$\Oz(\A_\ep(M_\ep, \B, \inv, \Theta))$ and $\Oz(\U_\ep(M_\ep, \B, \inv, \Theta))$ that is an alternative to that in \coref{Oz}. 

It uses the ideas of \cite{CGWZ2} for relating the ozone groups of prime PI algebras with the monoids of their normal elements. Recall that an element $u$ of a noncommutative ring $R$ is normal if $Ru =uR$. If $R$ is a domain (i.e., it has no zero divisors), then each of its nonzero normal elements $u$ induces an automorphism $\eta_u \in \Aut(R)$ such that 
\[
u r = \eta_u(r) u, \quad \forall \; r \in R.
\]
In other words $\eta_u(r) = u r u^{-1}$. Moreover, $\eta_u \in \Oz(R)$ because $\eta_u(z) = u z u^{-1} =z$ for all $z \in Z(R)$.

For every domain $R$ we have a homomorphism 
\begin{equation}
\label{eta}
\eta: N(R) \to \Oz(R) 
\end{equation}
where $N(R)$ is the monoid of nonzero normal elements of $R$. This homomorphism is rarely injective since it is trivial on $Z(R) \backslash \{ 0 \}$. We note that \cite{CGWZ2} uses the automorphism $r \mapsto u^{-1} r u$ instead, but that one does not define a homomorphism as in \eqref{eta}.

We will make use of the following key result of \cite{CGWZ2}:

\ble{CGWZ2} (Chan--Gaddis--Won--Zhang) \cite[Lemma 1.9]{CGWZ2} 
The homomorphism \eqref{eta} is surjective for every PI ring $R$ that is a domain. Furthermore, if $R$ is graded by an ordered abelian group, then the restriction of $\eta$ to the monoid of homogeneous nonzero normal elements is also surjective. 
\ele

Returning to the setting of a root of unity quantum cluster algebra $\A_\ep(M_\ep, \B)$, denote 
\[
\Nset^\fr := \bigoplus_{i \in \fr} \Nset e_i
\quad \mbox{and} \quad 
(\Zset/\ell)^\fr := \bigoplus_{i \in \fr} (\Zset/\ell) \ol{e_i}.
\]
Recall from Sect. \ref{2.2} that, for a root of unity toric frame $M_\ep$, $\Om_{M_\ep} : \Zset^N \times \Zset^N \to \Zset/\ell$ denotes both an additive bicharacter and its matrix in $M_N(\Zset/\ell)$. Using the latter convention, we set
\begin{equation}
\label{Om-fr}
\Im_{\fr} \big( \Om_{M_\ep} \big) := \Om_{M_\ep} \big( (\Zset/\ell)^\fr \big).
\end{equation}

The frozen variables $\{ M_\ep(e_i) \mid i \in \fr \}$ do not change under mutation; that is, $M'_\ep(e_i) = M_\ep(e_i)$ for all seeds $(M'_\ep, \B')$ that are obtained by mutations from the initial one $(M_\ep, \B)$. Hence, 
\[
M'_\ep(f) = M_\ep(f), \quad \forall f \; \in \Nset^\fr,
\]
and thus $M_\ep(f)$ for $ f \; \in \Nset^\fr$ is a normal element of each mutated quantum torus $\TT_\ep(M'_\ep)$ 
and the mixed quantum torus-affine subalgebra 
$\TT_\ep(M'_\ep)_{\geq}$
satisfying 
\begin{equation}
\label{Mepfg}
M_\ep(f) M'_\ep(g) = \ep^{\Om_{M'_\ep}(f,g)} M'_\ep(g) M_\ep(f), 
\quad \forall g \in \Zset^N.
\end{equation}
Therefore, it is a normal element of the algebras $\A_\ep(M_\ep, \B, \inv, \Theta)$ and $\U_\ep(M_\ep, \B, \inv, \Theta)$ for all collections of seeds $\Theta$. 

By the compatibility relation \eqref{compatible}, we have 
\[
\Im_{\fr} \big( \Om_{M_\ep} \big) \subseteq \Ker (\ol{\B}^\top). 
\]
So we can restrict the action from \prref{ep-act} to $\Im_{\fr} \big( \Om_{M_\ep} \big)$.

\ble{oz-normal} Assume that $\ell$ is a positive odd integer, $\ep^{1/2} \in \KK$ is a primitive $\ell$-th root of unity, 
$\Theta$ is a connected collection of seeds of the root of unity quantum cluster algebra $\A_\ep(M_\ep, \B)$, and $\inv \subseteq \fr$ is arbitrary. 

The restriction of the action of $\Ker(\ol{\B}^\top)$ on 
$\A_\ep(M_\ep, \B, \inv, \Theta)$ and 
$\U_\ep(M_\ep, \B, \inv, \Theta)$ from \prref{ep-act} to $\Im_{\fr} \big( \Om_{M_\ep} \big)$ defines injective group homomorphisms 
\[
\Im_{\fr} \big( \Om_{M_\ep} \big) \hra 
\Oz(\A_\ep(M_\ep, \B, \inv, \Theta)) 
\quad 
\mbox{and} \quad
\Im_{\fr} \big( \Om_{M_\ep} \big) \hra 
\Oz(\U_\ep(M_\ep, \B, \inv, \Theta)).
\]
\ele
\begin{proof}
Taking into account the skew-symmetry of 
$\Om_{M_\ep} : \Zset^N \times \Zset^N \to (\Zset/\ell)$, from \eqref{Mepfg} we obtain that 
\[
(\Om_{M_\ep} \ol{f} \cdot) = \eta_{M_\ep((\ell-1)f)}, \quad \forall \; f \in \Nset^\fr.
\]
These automorphisms of $\A_\ep(M_\ep, \B, \inv, \Theta)$ and 
$\U_\ep(M_\ep, \B, \inv, \Theta)$ fix their centers because of \eqref{eta}.
\end{proof}
\bth{normal-Oz} Assume that $\ell$ is a positive integer satisfying condition (Coprime), $\ep^{1/2} \in \KK$ is a primitive $\ell$-th root of unity, 
$\Theta$ is a connected collection of seeds of the root of unity quantum cluster algebra $\A_\ep(M_\ep, \B)$ that is a nerve, and $\inv \subseteq \fr$ is arbitrary. 

The injective homomorphisms from \leref{oz-normal} restrict to isomorphisms
\[
\Im_{\fr} \big( \Om_{M_\ep} \big) \cong 
\Oz(\A_\ep(M_\ep, \B, \inv, \Theta)) 
\quad 
\mbox{and} \quad
\Im_{\fr} \big( \Om_{M_\ep} \big) \cong 
\Oz(\U_\ep(M_\ep, \B, \inv, \Theta)).
\]

In other words, every element of the ozone groups of the two algebras is of the form 
$\eta_u$ for some $u$ that is a product of the frozen variables $M_\ep(e_i)$, $i \in \fr$; all powers can be taken in the range $[0,\ell-1]$ since the $\ell$-th powers of all frozen variables are central. 
\eth
\begin{proof} We show that the homomorphisms in \leref{oz-normal} are surjective. Let $R:= \A_\ep(M_\ep, \B, \inv, \Theta)$ and $R:= \U_\ep(M_\ep, \B, \inv, \Theta)$ in the two cases, respectively. 

Fix $\phi \in \Oz(R)$. In both cases $\TT_\ep(M_\ep)$ is a central localization of $R$:
\[
R[M_\ep(\ell e_k)^{-1} : k \in \ex \sqcup \ninv] = \TT_\ep(M_\ep).
\]
Therefore, $\phi$ extends to an element of the ozone group of this localization. The root of unity quantum torus $\TT_\ep(M_\ep)$ is $\Zset^N$-graded by 
\[
\deg M_\ep(f) :=f \quad \forall f \in \Zset^N
\]
and its homogeneous components are all one-dimensional. 
Applying the second statement in \leref{CGWZ2} to it gives that there exists $f \in \Zset^N$ such that
\begin{equation}
\label{rf}
\phi(r) = M_\ep(-f) r M_\ep(f), \quad \forall \; r \in R. 
\end{equation}
Hence, $\phi( M_\ep (e_i)) = \ep^m M_\ep(e_i)$ for some $m \in \Zset/\ell$ depending on $i$. 
The same can be said about every seed $(M'_\ep, \B') \in \Theta$.
Therefore, $\phi$ satisfies condition \eqref{condII} in \thref{nerve}. By \thref{nerve}, there exists $\theta \in \Ker \big( \ol{\B}^\top \big)$ such that 
\[
\phi ( M_\ep (e_i) ) = \ep^{(\theta)_i} M_\ep(e_i), \quad \forall \; i \in [1,N].
\]
However, from \eqref{rf}, we have 
\[
\phi ( M_\ep(e_i)) = \ep^{(\Om_{M_\ep} \ol{f})_i} M_\ep(e_i), \quad
\forall \; i \in [1,N].
\]
Therefore 
\[
\Om_{M_\ep} \big( \ol{f} \big) = \theta \in \Ker \big( \ol{\B}^\top \big)
\quad \mbox{and} \quad 
\phi = \big( \Om_{M_\ep}(\ol{f}) \cdot \big).
\]
Eq. \eqref{compatible} implies
\[
\ol{\B}^\top \big( \Om_{M_\ep}(\ol{f}) \big) = \sum_{k \in \ex}
\ol{d_k} (\ol{f})_k \ol{e_k},
\]
where $\{ d_k \mid k \in \ex\}$ are the skew-symmetrizing integers for the principal part of the exchange matrix $\B$. 
Condition (Coprime) gives that $\ol{d_k}$ are units of $\Zset/\ell$ for all $k \in \ex$. So, $(\ol{f})_k =0$ for all $k \in \ex$. Therefore, $\ol{f} \in (\Zset/\ell)^\fr$. Since $\phi = \big( \Om_{M_\ep}(\ol{f}) \cdot \big)$, $\phi$ is in the image of the injective homomorphism in \leref{oz-normal}, which completes the proof of the theorem.
\end{proof}

\bex{non-cycl}
    Let $\ell=9$ and $\KK$ be a field with $\charr \KK \neq 3$ with $\ep^{1/2}$ a primitive $9$th root of unity. Take $N=4$, $\ex = \{1\}$, $\fr = \{2,3,4\}$ and let 
    \begin{equation}
    \Omega :=\begin{pmatrix}
             0 & 1 & 0 & 0\\
             -1 & 0 & 0 & 0\\
             0 & 0 & 0 & 3\\
             0 & 0 & -3 &0
            \end{pmatrix} \in M_4(\Zset/9), \quad \widetilde{B} := \begin{pmatrix}
                0\\-1\\0\\0
            \end{pmatrix} \in M_{4\times \ex}(\Zset).
    \end{equation}
    Let $\FF :=\Fract(\TT_\ep(\Omega))$ and $M_\ep : \Zset^4 \rightarrow \FF$, $M_\ep(\nu) := x^\nu$ be the toric frame such that $\Omega_{M_\ep} = \Omega$. Then $(M_\ep, \widetilde{B})$ is a root of unity quantum seed, condition (Coprime) holds, and $$\Oz(\A_\ep(M_\ep,\widetilde{B},\inv,\Theta)) = \Oz(\U_\ep(M_\ep,\widetilde{B},\inv,\Theta)) \cong \Zset/9 \oplus \Zset / 3 \oplus \Zset / 3$$
    for every $\inv \subseteq \fr$ and every nerve $\Theta$. In particular, the ozone group is not a direct sum of cyclic groups of the same order, while by \thref{torus} the corresponding Frobenius automorphism group is isomorphic to $(\Zset/9)^3$.
\eex

Comparing \coref{Oz}(ii) and \thref{normal-Oz} gives the following general result for the properties of exchange matrices of cluster algebras that is of independent interest:
\bco{equality} Assume that $\KK$, $\ell$, $\ep^{1/2}$, $\B$, $\Theta$ and $\inv$ are as in \thref{normal-Oz}.
Then we have the equality 
\[
\Ker \big( \, \ol{\B}^\top \big) \cap \big( \ol{\Ker \Om_{M_\ep} }\big)^\perp = \Im_{\fr} \big( \Om_{M_\ep} \big).
\]
\eco

\subsection{Type $A_2$ with arbitrary coefficients}
In this subsection we determine the ozone groups of all root of unity quantum cluster algebras whose principal part is of type $A_2$. Set
\[
\ex := \{1,2\}, \qquad \fr := [3,N], \qquad m := N-2,
\qquad
B := \begin{pmatrix} 0 & 1\\ -1 & 0 \end{pmatrix}.
\]
Consider the exchange matrices
\begin{equation}
\label{eA2B}
\widetilde{B} \; = \; \begin{pmatrix} B \\ c \end{pmatrix}
\in M_{N \times \ex}(\Zset),
\qquad c \in M_{m \times 2}(\Zset) \ \text{arbitrary},
\end{equation}
with principal part the skew-symmetric matrix $B$ of type $A_2$.
Since $D = I_2$, condition (Coprime) is satisfied exactly when $\ell$ is odd and coprime to $\charr \KK$. We write the bilinear form of a root of unity toric frame $M_\ep$ in the
block form associated to the decomposition $[1,N] = \ex \sqcup \fr$,
\begin{equation}
\label{eA2Om}
\Omega_{M_\ep} \; = \;
\begin{pmatrix} \Lambda & Q \\ -Q^{\top} & S \end{pmatrix},
\qquad
\Lambda \in M_2(\Zset/\ell), \quad
Q \in M_{2 \times m}(\Zset/\ell), \quad
S \in M_m(\Zset/\ell),
\end{equation}
where the matrices $\Lambda$ and $S$ are skew-symmetric.
\ble{A2compat}
Let $\widetilde{B}$ be as in \eqref{eA2B} and let $\Omega_{M_\ep}$ be as
in \eqref{eA2Om}. The pair $(\Omega_{M_\ep}, \widetilde{B})$ is
$\ell$-compatible (with respect to $D = I_2$) if and only if
\begin{equation}
\label{eA2compat}
Q = - B \, \overline{c}^{\,\top} S
\qquad\text{and}\qquad
\Lambda = (1+\sigma) B, \quad \text{where} \quad
\sigma := \bigl( \overline{c}^{\,\top} S \, \overline{c} \bigr)_{12} \in \Zset/\ell .
\end{equation}
\ele
\begin{proof}
    $\ol{\widetilde{B}}^{\top} \Omega_{M_\ep} = \begin{pmatrix} I_2 & 0\end{pmatrix}$
is equivalent to the pair of identities
\begin{equation}
\label{eA2two}
B^{\top} \Lambda - \overline{c}^{\,\top} Q^{\top} = I_2
\qquad\text{and}\qquad
B^{\top} Q + \overline{c}^{\,\top} S = 0 .
\end{equation}
Since $B^2 = -I_2$, we have $B^{\top} = -B = B^{-1}$, so the second identity holds if and only if
$$Q = - (B^{\top})^{-1} \overline{c}^{\,\top} S = - B \overline{c}^{\,\top} S.$$ Since $S$ is skew-symmetric, we get
$Q^{\top} = - S^{\top} \overline{c} \, B^{\top} = - S \overline{c} \, B$, and thus, 
\[
\overline{c}^{\,\top} Q^{\top}
= - \bigl( \overline{c}^{\,\top} S \, \overline{c} \bigr) B
= - \sigma B^2 = \sigma I_2 ,
\]
where in the second equality we used that the matrix $\overline{c}^{\,\top} S \, \overline{c}$ is skew-symmetric, and thus equals $\sigma B$ for $\sigma := (\overline{c}^{\,\top} S \, \overline{c})_{12}$. Therefore, the first identity in \eqref{eA2two} is equivalent to
$B^{\top} \Lambda = (1+\sigma) I_2$, which is equivalent to
$\Lambda = (B^{\top})^{-1} (1+\sigma) I_2 = (1+\sigma) B$.
\end{proof}

\bpr{A2}
Assume that $\ell$ is odd and coprime to $\charr \KK$ and let $\widetilde{B}$ be an exchange matrix of type
$A_2$ with entries as in \eqref{eA2B}. For a root of unity quantum seed 
$(M_\ep, \widetilde{B})$ write $\Omega_{M_\ep}$ in the block
form \eqref{eA2Om}. Then for every $\inv \subseteq \fr$ and every nerve $\Theta$,
$$
\Oz\bigl( \A_\ep(M_\ep, \widetilde{B}, \inv, \Theta) \bigr)
=\Oz\bigl( \U_\ep(M_\ep, \widetilde{B}, \inv, \Theta) \bigr)\cong  \Im(S).$$
\epr
\begin{proof}
    By Theorem~\ref{tnormal-Oz}, both ozone groups equal
$\Im_{\fr}(\Omega_{M_\ep})$. By \eqref{eA2Om},
$$\Im_{\fr}(\Omega_{M_\ep})
= \bigl\{ (Qg, Sg) \mid g \in (\Zset/\ell)^m \bigr\} $$
and by Lemma~\ref{lA2compat}, $Qg = - B \overline{c}^{\,\top} (S g)$, hence the map
$\Im(S) \to \Im_{\fr}(\Omega_{M_\ep})$, given by
$v \mapsto (- B \overline{c}^{\,\top} v, v)$, is a group isomorphism.
\end{proof}
\bth{A2realize}
Let $H$ be a finite abelian group whose order is odd and invertible in $\KK$. Write
\[
H \; \cong \; \Zset/a_1 \oplus \cdots \oplus \Zset/a_k,
\qquad a_1 \mid a_2 \mid \cdots \mid a_k,
\]
set $\ell := a_k$, and assume that $\ep^{1/2}$ is a primitive $\ell$-th root of unity. Denote $N := 2k+2$, $m := 2k$. Choose an arbitrary $c \in M_{m \times 2}(\Zset)$ 
and set
\begin{equation}
\label{eSrealize}
\widetilde{B} := \begin{pmatrix} B \\ c \end{pmatrix}
\in M_{N \times \mathrm{\ex}}(\Zset),
\qquad
S := \bigoplus_{j=1}^{k} \frac{\ell}{a_j}\, B \; \in \; M_{m}(\Zset/\ell).
\end{equation}
Let $(M_\ep, \widetilde{B})$ be the root of unity quantum seed attached to the
pair $(c,S)$ by Lemma~\ref{lA2compat}. Then $\ell$ satisfies condition
\textup{(Coprime)}, and for every $\inv \subseteq \fr$ and every nerve $\Theta$,
\[
\Oz \bigl( \A_\ep(M_\ep, \widetilde{B}, \inv, \Theta) \bigr)
\; = \;
\Oz \bigl( \U_\ep(M_\ep, \widetilde{B}, \inv, \Theta) \bigr)
\; \cong \; H \oplus H.
\]
\eth

\begin{proof}
Since $\ell = \exp(H)$ divides $|H|$ and the latter is odd and invertible
in $\KK$ by hypothesis, we have that (Coprime) holds. By \leref{A2compat} we can construct an $\ell$-compatible pair $(\Omega_{M_\ep}, \widetilde{B})$ with frozen block skew-symmetric matrix $S$, which is a direct sum of scalar multiples of $B$.  
Write $s_j := \ell/a_j$. Since $S$ is block diagonal
and $B$ is invertible over $\Zset/\ell$, we have
\[
\Im(S)  =  \bigoplus_{j=1}^{k} \Im \bigl( s_j B \bigr)
 = \bigoplus_{j=1}^{k} \bigl( s_j \Zset/\ell \bigr)^{2}
 \cong  \bigoplus_{j=1}^{k} \bigl( \Zset/(\ell/s_j) \bigr)^{2}
 = \bigoplus_{j=1}^{k} \bigl( \Zset/a_j \bigr)^{2}
\cong H \oplus H .
\]
\prref{A2} now gives the desired result.
\end{proof}
\bex{ozone39}
Let $H = \Zset/3 \oplus \Zset/9$, so $\ell = 9$, $k = 2$, $N = 6$. Taking $c=0$
and $S = 3B \oplus B$ gives
\[
\Oz \bigl( \A_\ep(M_\ep, \widetilde{B}, \inv, \Theta) \bigr)
\; \cong \; (\Zset/3)^2 \oplus (\Zset/9)^2 .
\]
\eex

\sectionnew{Frobenius and ozone automorphism groups of quantum groups}
\label{qgroups}
In this section we use the general theorems from the previous sections to classify the Frobenius and 
ozone automorphism groups of the quantized function
algebras of complex simple algebraic groups at roots of unity. 
\subsection{Integral forms of quantum groups and specializations at roots of unity}
\label{qg-forms}
Let $\g$ be a finite dimensional, complex, simple Lie algebra with Cartan matrix $C=(c_{ij})_{i,j=1}^r$ and $G$ be the corresponding connected, simply connected, complex algebraic group. 
Denote by $P$ the weight lattice of $\g$, by $P_+$ the set of its dominant integral weights, by $\{\vpi_i \mid i \in [1,r] \}$ its fundamental weights, and by $W$ its Weyl group (with length function $l : W \to \Nset$).

Let $U_q(\g)$ be the Drinfeld--Jimbo quantized enveloping algebra of $\g$ over $\Cset(q)$; we refer the reader to \cite{BrG,CP} for details. For $\la \in P_+$, denote by $L(\la)$ its (finite dimensional) highest weight module of highest weight $\la$, and for $\mu \in P$, its corresponding weight space $L(\la)_\mu$.
Decompose $L(\la) = \oplus_{\mu \in P} L(\la)_\mu$ and identify $L(\la)^* \cong \oplus_{\mu \in P} (L(\la)_\mu)^*$.
Fix a highest weight vector $v_\la \in L(\la)_\la$ and extremal weight vectors $v_{w \la} \in L(\la)_{w \la}$ for $w \in W$ obtained from $v_\la$ by Lusztig's braid group action. Recall that $\dim L(\la)_{w \la}=1$ and fix $\xi_{w \la} \in (L(\la)_{w \la})^*$ such that $\langle \xi_{w \la}, v_{w \la} \rangle =1$. The quantized coordinate ring of $G$ is the Hopf subalgebra $\OO_q(G) \subseteq (U_q(\g))^*$
spanned by the matrix coefficients of the modules $L(\la)$ for $\la \in P_+$. As a vector space it is given by the quantum Peter--Weyl decomposition 
\[
\OO_q(G) \cong \bigoplus_{\la \in P_+} L(\la)^* \otimes L(\la).
\]
The summand for $\la \in P_+$ is spanned by the matrix coefficients $c_{\xi, v}^\la$ corresponding to $\xi \in L(\la)^*$ 
and $v \in L(\la)$. The algebra $\OO_q(G)$ is 
$P \times P$-graded by 
\begin{equation}
\label{grading}
\OO_q(G)_{\mu, \nu} = \bigoplus_{\la \in P_+, \, \xi \in (L(\la)_\mu)^*, \, v \in L(\la)_\nu} c_{\xi, v}^\la. 
\end{equation}

Denote 
\[
\cA := \Cset[q^{\pm 1}] 
\quad \mbox{and} \quad 
\cAt := \Cset[q^{\pm 1/2}].
\]
De Concini, Kac and Procesi \cite{DK,DKP} defined an $\cA$-form $U_q(\g)_{\cA}$ of $U_q(\g)$, called the {\em{big quantum group}}. De Concini and Lyubashenko \cite{DCL94} constructed the corresponding $\cA$-form $\OO_q(G)_{\cA}$ of the quantized coordinate ring $\OO_q(G)$. In \cite{DK,DCL94} integral forms were defined over $\Zset[q^{\pm 1}]$; we will only need the forms over $\cA$. The generalized quantum minors of $G$ are the matrix coefficients 
\[
\Delta_{u \vpi_i, w \vpi_i} := c_{\xi_{u \vpi_i}, v_{w \vpi_i}}^{\vpi_i}
\]
for $i \in [1,r]$ and $w, u \in W$. 
They belong to the $\cA$-form $\OO_q(G)_{\cA}$.
At $q=1$, they specialize to the generalized minors $\Delta'_{u \vpi_i, w \vpi_i} \in \OO(G)$. More generally, set 
\begin{equation}
\label{la}
\Delta_{u \la, w \la} := c_{\xi_{u \la}, v_{w\la}}^{\la}
\end{equation}
for $\la \in P_+$ and $w, u \in W$. 

The quantized coordinate ring of $G$ at a primitive $\ell$-th root of unity $\ep$ is defined \cite{DCL94} as the $\Cset$-algebra obtained by the specialization 
\begin{equation}
\label{OepG}
\OO_\ep(G) := \OO_q(G)_{\cA}/ (q - \ep)
\cong \OO_q(G)_{\cA} \otimes_{\cA} \Cset.
\end{equation}
Here $(q-\ep)$ denotes the centrally generated ideal of $\OO_q(G)_{\cA}$ by $q - \ep$ and the tensor product is taken with respect to the $\Cset$-linear homomorphism $\cA \to \Cset$ sending $q$ to $\ep$. 

By abuse of notation, we will denote the image of an element $y \in \OO_q(G)$ in $\OO_\ep(G)$ by the same symbol. 

De Concini and Lyubashenko constructed a quantum Frobenius map 
\[
\Fr_\ep : \OO(G) \hookrightarrow Z(\OO_\ep(G))
\]
as the dual of Lusztig's quantum Frobenius map \cite{Lus2}, a surjective homomorphism from Lusztig's specialization of the divided power integral form \cite{Lus} of $U_q(\g)$ to $U(\g)$. Denote the image of the quantum Frobenius map by $\cC_\ep(G)$:
\begin{equation}
\label{qFr}
\Fr_\ep : \OO(G) \stackrel{\cong}{\longrightarrow} \cC_\ep(G) \subset Z( \OO_\ep(G)).
\end{equation}

\subsection{Cluster structures on quantum groups at roots of unity}
\label{qg-cluster}
Denote by $w_0$ the longest element of the Weyl group $W$ and set $n:= \dim \g - r$. Consider a reduced word $\ubi=(i_1,\ldots,i_n)$ for $(w_0,w_0) \in W \times W$ such that the simple reflections for the first component are labeled by $[-r,-1]$ and those for the second component by $[1,r]$. (Recall that a reduced word for an element of a Coxeter group is the list of the indices of a reduced expression of the element in terms of simple reflections.) 
Denote by $(u_{\leq k},w_{\leq k})\in W\times W$ the element corresponding to the product of the simple reflections indexed by $(i_1,\ldots,i_k)$ in this order. Let $(i_{-r},\ldots,i_{-1})$ be a permutation of $[1,r]$. 

The Berenstein--Fomin--Zelevinsky seeds \cite{BFZ05}
for $\OO(G)$ are indexed by reduced expressions $\ubi$
of $(w_0, w_0)$. The set of their exchangeable and frozen indices is chosen to be $I=[-r,-1]\sqcup [1,n]$.
This slightly differs from the conventions in Sect. \ref{prelim} where we used intervals $[1,N]$ for such indexing, but is the standard choice in this setting.
The exchange matrix $\B(\ubi)$ is defined as follows. For $k \in I$, set
$k[1]:=\min(\{j\in I|j>k,|i_j|=|i_k| \}\sqcup \{\infty\})$. The set of unfrozen indices is chosen to be $\ex(\ubi):=\{k\in [1,n]|k[1]\leq n\}$ and the frozen ones $\fr(\ubi):= I \backslash \ex(\ubi)$; mutations are considered up to relabeling of indices and this is why $\ex$ and $\fr$ depend on $\ubi$. Denote $s_k:=\sign(i_k)$. The entries of the exchange matrix $\B(\ubi):=(b_{jk}(\ubi))_{j\in I, k\in \ex(\ubi)}$ are chosen to be
\begin{align*}
    b_{jk}(\ubi)=
    \begin{cases}
        -s_k, & k=j[1] \\
        s_j, & j=k[1]\\
        -s_k c_{|i_j|,|i_k|}, & s_k=s_{j[1]},\quad j<k<j[1]<k[1]\\
        s_j c_{|i_j|,|i_k|}, & s_j=s_{k[1]},\quad k<j<k[1]<j[1]\\
        -s_k c_{|i_j|,|i_k|}, & s_k=-s_{k[1]},\quad j<k<k[1]<j[1]\\
        s_j c_{|i_j|,|i_k|}, & s_j=-s_{j[1]},\quad k<j<j[1]<k[1]\\
        0, & \text{otherwise.}       
    \end{cases}
\end{align*}
For $k\in I$, let
\begin{align*}
    \gamma_k:=\begin{cases}
        \varpi_{i_k},&\text{ if }k\in [-r,-1]\\
        u_{\leq k}\varpi_{|i_k|},&\text{ if }k\in [1,n]
    \end{cases},\  \delta_k:=\begin{cases}
        w_0^{-1}\varpi_{i_k},&\text{ if }k\in [-r,-1]\\
        w_0^{-1}w_{\leq k}\varpi_{|i_k|},&\text{ if }k\in [1,n]
    \end{cases}
\end{align*}
and 
\[
x(\ubi)_k:= \Delta'_{\gamma_k, \delta_k}, \quad k \in I.
\]
The Berenstein--Fomin--Zelevinsky \cite{BFZ05} seeds of $\OO(G)$ are the seeds 
\[
\Sigma(\ubi):= (\x(\ubi), \B(\ubi)), \quad \mbox{where} \quad \x(\ubi):= \{ x(\ubi)_k \mid k \in I \}. 
\]
Throughout the section we choose $\inv = \varnothing$. 

By \cite[Theorem B.1]{QY} and \cite[Theorem 3.1]{O}, respectively, we have 
\begin{equation}
\label{OG}
\OO(G) = \U(\B, \varnothing) \; \mbox{for all types $G$}, \; \;    
\OO(G) = \A(\B, \varnothing) = \U(\B, \varnothing)  \; \mbox{for $G$ not of type $F_4$}.
\end{equation} 

Berenstein and Zelevinsky \cite{BerZe} conjectured that for all reduced words $\ubi$ of $(w_0, w_0)$, there is a (unique) toric frame of $\OO_q(G)_{\cAt} := \OO_q(G) \otimes_{\cA} \cAt$ 
with
\[
M_q(\ubi)(e_k):= \Delta_{\gamma_k, \delta_k}, \quad \forall \; k \in I 
\]
and with matrix $\Lambda(\ubi):= (\lambda(\ubi)_{kj})_{k,j \in I}$ such that $\lambda(\ubi)_{kj}:=(\gamma_k,\gamma_j)-(\delta_k,\delta_j),\ \forall k>j$. 
This was proved in \cite{GY1,QY}. Furthermore, by \cite[Theorem A]{OQY26}, 
\begin{equation}
\label{OqAU}
\OO_q(G)_{\mathbb{A}} = \U_q(M_q(\ubi), \B(\ubi), \varnothing)_{\mathbb{A}} =  \A_q(M_q(\ubi), \B(\ubi), \varnothing)_{\mathbb{A}}
\end{equation}
if $G$ is not of type $F_4$ or $E_8$, where  
\[
{\mathbb{A}} := 
\begin{cases}
\cAt[(q^2 +1)^{-1}], &\mbox{if $G$ is of type $G_2$}
\\
\cAt, & \mbox{otherwise}
\end{cases}
\]
and  $\OO_q(G)_{\mathbb{A}} := \OO_q(G) \otimes_{\cA} {\mathbb{A}}$. 
In \eqref{OqAU}, the second and third terms denote the upper quantum cluster and quantum cluster algebras in the sense of Berenstein--Zelevinsky 
\cite{BerZe} defined over ${\mathbb{A}}$. The first equality holds for all types over $\cAt$, \cite[Theorem A]{OQY26}. 

In \cite{GY1,QY} it was proved that for every reduced word $\ubi$ of $(w_0,w_0)$, 
the cluster variables $M_q(\ubi)(e_i)$ generate a (left and right) denominator set, and the corresponding localization 
of $\OO_q(G)$ is the quantum torus generated by them. The specialization \eqref{OepG} gives that the same is true for their images in $\OO_\ep(G)$. 
Therefore, for all reduced words $\ubi$ of $(w_0, w_0)$, $\OO_\ep(G)$ has a root of unity toric frame $M_\ep(\ubi)$ 
with 
\[
M_\ep(\ubi)(e_k):= \Delta_{\gamma_k, \delta_k}, \quad \forall \; k \in I 
\]
and with matrix $\ol{\Lambda(\ubi)}$. It gives rise to the root of unity quantum seed $\Sigma_\ep:=(M_\ep(\ubi), \B(\ubi))$ of $\OO_\ep(G)$. 

\bth{qCl} Let $G$ be a connected, simply connected, complex simple algebraic group not of type $F_4$ or $E_8$, $\ell$ be an odd positive integer
not divisible by 3 if $G$ is of type $G_2$, and $\ep^{1/2}$ be a primitive $\ell$-th root of unity. Then the following hold:
\begin{enumerate}
\item[(i)] For all reduced words $\ubi$ of $(w_0, w_0)$ the root of unity quantum seeds $\Sigma_\ep:=(M_\ep(\ubi), \B(\ubi))$ of $\OO_\ep(G)$ 
are obtained from each other by mutations and 
\begin{equation}
\label{OepA}
\OO_\ep(G) =  \A_\ep(M_\ep(\ubi), \B(\ubi), \varnothing).
\end{equation}
\item[(ii)] The Frobenius centers of the algebras in the first part coincide: 
\[
\cC_\ep(G) =  \C\A_\ep(M_\ep(\ubi), \B(\ubi), \varnothing).
\]
\end{enumerate} 
\eth
\begin{proof} The assumption on $\ell$ implies that condition (Coprime) is satisfied 
because the skew-symmetrizing integers of the principal part of $\B(\ubi)$ are chosen 
to be $d_{|i_j|}$ where $\{d_i \mid i \in [1,r] \}$ are the symmetrizing (positive, relatively prime) 
integers for the Cartan matrix of $G$. 

(i) The specialization of quantum cluster algebras \cite[Eq. (5.2)]{NTY2} maps 
quantum mutations to root of unity quantum mutations.
This and the fact that the quantum seeds $(M_\ep(\ubi), \B(\ubi))$ are obtained from each other by mutations \cite[Theorem A]{QY} 
imply the first statement in part (i). Eq. \eqref{OepA} follows from applying the specialization result \cite[Theorem 5.11]{NTY2} and \eqref{OqAU}. 

(ii) The Frobenius map \eqref{qFr} satisfies 
\[
\Fr_\ep( \Delta'_{u \vpi_i, w \vpi_i} ) = \Delta_{u \vpi_i, w \vpi_i}^\ell.
\]
Under condition (Coprime) the $\ell$-th powers of root of unity quantum cluster variables mutate at the corresponding cluster variables by \cite[Proposition 4.4]{NTY2}.
Part (ii) follows from these two facts and the equalities in \eqref{OG} and \eqref{OepA}. 
\end{proof}

\subsection{The full center of $\OO_\ep(G)$}
In the appendix to \cite{DCL94} Enriquez proved that the full center of $\OO_\ep(G)$ is given by 
\begin{equation}
\label{ZO}
Z(\OO_\ep(G)) \cong \cC_\ep(G) \otimes_{\cC_\ep(G) \cap Z_q} Z_q,
\end{equation}
where $Z_q$ is the central subalgebra of $\OO_q(G)_\cA$ generated by the binomials
\begin{equation}
\label{Zq}
Z_q
=
\Cset \left[
\left(\Delta_{w_0\varpi_i,\varpi_i}\right)^k
\left(\Delta_{-w_0\varpi_i,-\varpi_i}\right)^{\ell-k}
\; : \; 
i\in[1,r],\ 0\leq k\leq \ell
\right].
\end{equation}
The $2r$ quantum minors that appear in these binomials are precisely the frozen variables in the root of unity Berenstein--Zelevinsky seeds $\Sigma_\ep(\ubi)$.  
\subsection{Relative automorphism groups and ozone groups of quantum function algebras at roots of unity}
Let $T:= \Gset_m^r$. The torus $T^2$ acts on $\OO_\ep(G)$ by algebra automorphisms by 
\begin{equation}
\label{TactO}
(t_1, \ldots, t_r, t_{r+1}, \ldots, t_{2r}) \cdot c = t_1^{m_1} \cdots t_r^{m_r} t_{r+1}^{n_1} \cdots t_{2r}^{n_r} c
\end{equation}
for
\[
c \in \OO_\ep(G)_{\mu, \nu}, \quad - w_0(\mu) = m_1 \vpi_1 + \cdots + m_r \vpi_r, \; \; 
\nu = n_1 \vpi_1 + \cdots + n_r \vpi_r, \; m_i, n_i \in \Zset,
\]
recall the $P \times P$-grading \eqref{grading} of $\OO_\ep(G)$. Denote the diagonal embedding
\[
\Delta : T \hra T^2.
\]
We will use the same notation for the map on the $\ell$-torsion subgroup of its $\KK$-points: $\Delta : C_\ell^r \hra C_\ell^{2 r}$.

The main classification result in this section is the following:

\bth{Oepgroup}
Let $G$ be a connected, simply connected, complex simple algebraic group of rank $r$ not of type $F_4$ or $E_8$, $\ell$ be an odd positive integer
not divisible by 3 if $G$ is of type $G_2$, and $\ep$ be a primitive $\ell$-th root of unity. 

\begin{enumerate}
\item[(i)] The Frobenius automorphism group of $\OO_\ep(G)$  is given by 
\[
\Aut(\OO_\ep(G)/\cC_\ep(G)) = T^2(\KK)[\ell] = C_\ell^{2r}.
\]

Every endomorphism of $\OO_\ep(G)$ fixing $\cC_\ep(G)$ is an automorphism.

\item[(ii)] The ozone group of $\OO_\ep(G)$ is
\[
\Oz(\OO_\ep(G)) = \Delta(T)(\KK)[\ell] =  \Delta(C_\ell^{r}) \cong C_\ell^r.
\]
The elements of the group $\Oz(\OO_\ep(G))$ are precisely the automorphisms of $\OO_\ep(G)$ of the form 
\[
r \mt 
\Delta_{w_0 \la, \la} r \Delta_{w_0 \la, \la}^{-1} \quad \mbox{for} \quad\la = n_1 \varpi_1 + \cdots + n_r \varpi_r, \; \mbox{with} \; n_i \in [0, \ell -1],
\]
recall \eqref{la}.
\end{enumerate}
\eth
The equalities in part (i) and the first statement in part (ii) mean that the elements of the two automorphism groups come from the action of the $\ell$-torsion subgroups of $T^2$ and $\Delta(T)$.

\begin{proof} (i) Due to the two equalities in \thref{qCl}, we can match the relative automorphism 
groups of the setting of the quantum function algebras at roots of unity and the setting of root of unity quantum cluster algebras:
\[
\Aut(\OO_\ep(G)/ \cC_\ep(G)) = 
\Aut ( \A_\ep(M_\ep(\ubi), \B(\ubi), \varnothing) /  \C \A_\ep(M_\ep(\ubi), \B(\ubi), \varnothing)).
\]
The character lattice of the $2r$-torus $T^2$ is isomorphic to $P \times P$ and this isomorphism is compatible with the setup for its action \eqref{TactO} on $\OO_\ep(G)$. 
From this we obtain that the action \eqref{TactO} is faithful. Moreover, the rank of the torus $T^2$ equals $2r$, which is the number of frozen variables of the 
root of unity quantum cluster algebra $\A_\ep(M_\ep(\ubi), \B(\ubi))$. Hence, we can apply \thref{torus}(i), which implies the first statement in part (i). 
The second statement follows from \thref{aut-qcl}(i).  

(ii) \thref{torus}(i) and eqs. \eqref{ZO} and \eqref{Zq} imply that $\Oz(\OO_\ep(G))$ consists of those $t:=(t_1, \ldots, t_r, t_{r+1}, \ldots, t_{2r}) \in C_\ell^{2r}$ such that
$t \cdot y_{i,k} = y_{i,k}$ for all $i \in [1,r]$, $k \in [0, \ell]$, where 
\[
y_{i,k} := 
 \left(\Delta_{w_0\varpi_i,\varpi_i}\right)^k \left(\Delta_{-w_0\varpi_i,-\varpi_i}\right)^{\ell-k}.
\]
Taking into account that $t_i^\ell = 1$, from \eqref{TactO}, we obtain 
\[
t \cdot y_{i,k} = t_i^{-k} t_{r+i}^k  y_{i,k}.
\]
Therefore, the condition $t \cdot y_{i,k} = y_{i,k}$ for all $i \in [1,r]$ is equivalent to $t_i = t_{r+i}$ for all $i \in [1,r]$, i.e., 
$t \in \Delta(T(\KK)[\ell]) = \Delta(C_\ell^r)$. 

The frozen variables of the toric frames $M_q(\ubi)$, and thus of the root of unity toric frames $M_\ep(\ubi)$, are
\begin{equation}
\label{fr-var}
\{ \Delta_{w_0\varpi_i,\varpi_i}, \Delta_{-w_0\varpi_i,-\varpi_i}
\mid i \in [1,r] \}. 
\end{equation}
We have \cite[Eqs. (2.22)--(2.23)]{Y-spec}
\begin{align*}
\Delta_{w_0\varpi_i,\varpi_i} c &= 
\ep^{\langle w_0 \varpi_i, \mu \rangle - 
\langle \varpi_i, \nu \rangle} 
c \Delta_{w_0\varpi_i,\varpi_i},
\\
\Delta_{-w_0\varpi_i,-\varpi_i} c & = 
\ep^{\langle w_0 \varpi_i, \mu \rangle -
\langle \varpi_i, \nu \rangle} 
c\Delta_{-w_0\varpi_i,-\varpi_i} 
\end{align*}
for all $c \in \OO_q(G)_{\mu, \nu}$. Recall \eqref{eta}. The above commutation relations show that the automorphisms $\eta_u$ in the ozone group of $\OO_\ep(G)$ associated to the monomials in the second set of frozen variables in \eqref{fr-var} are the same as the automorphisms $\eta_u$ in the ozone group of $\OO_\ep(G)$ associated to the monomials in the first set of frozen variables. \thref{normal-Oz} implies that the elements of $\Oz(\OO_\ep(G))$ are precisely the elements of the form $\eta_u$ for the normal elements 
\begin{equation}
\label{la2}
u = \Delta_{w_0 \varpi_1, \varpi_1}^{n_1} \cdots 
\Delta_{w_0 \varpi_r, \varpi_r}^{n_r}
\end{equation}
for $n_1, \ldots, n_r \in [0, \ell -1]$. The element in \eqref{la2} is a scalar multiple of 
\[
\Delta_{w_0 \la, \la} \quad \mbox{for} \quad\la = n_1 \varpi_1 + \cdots + n_r \varpi_r
\]
by \cite[Eq. (2.18)]{Y-spec}.
\end{proof} 
\sectionnew{Frobenius and ozone automorphism groups of quantum unipotent cells}
\lb{unipotent}
The quantum function algebras $\OO_\ep(G)$ of Sect. \ref{qgroups} are built from finite dimensional simple Lie algebras. In this section, we extend to the case of symmetrizable Kac--Moody algebras and classify the Frobenius and ozone automorphism groups of all quantum unipotent cell algebras at roots of unity $A_\ep(\n_+(w))$ associated to an arbitrary symmetrizable Kac--Moody algebra $\g$ and a Weyl group element $w$. 
\subsection{Symmetrizable Kac--Moody algebras and their quantized enveloping algebras}
\label{unip-KM}
We extend the setting from the last section to include symmetrizable Kac--Moody algebras $\g$, retaining the notations $(c_{ij})_{i,j=1}^r$,
$\{d_i\}_{i=1}^r$, $Q$, $P$, $P_+$, $W$, 
$\{\al_i\}_{i=1}^r$, 
$\{\vpi_i\}_{i=1}^r$ 
for its Cartan matrix, symmetrizing integers, root and weight lattices, dominant integral weights, Weyl group, simple roots and fundamental weights. Denote by $P^\vee$ and $\{ h_i\}_{i=1}^r$ its coweight lattice and set of simple coroots. There is a nondegenerate symmetric $\Qset$-valued bilinear form $(\cdot, \cdot)$ on $P$ such that
\[
\langle \mu, h_i \rangle = \frac{2 (\al_i, \mu)}{(\al_i, \al_i)} \quad \mbox{and} \quad (\al_i, \al_i) = 2 d_i \quad \mbox{for all} \; i \in [1,r], \mu \in P.
\]
We will use the same notation for the induced form on the coweight lattice $P^\vee$. 

Let $U_q(\g)$ be the quantized universal enveloping algebra of $\g$ over $\Cset(q)$, generated by $e_i$, $f_i$ and $q^h$ for $i \in [1,r]$, $h \in P^\vee$; we refer the reader to \cite[Sect. 1.1]{Kas91} for full details. Denote $q_i : = q^{d_i}$ and let $[k]_{q_i}$ be the $q$-integer with base $q_i$. 

The unital subalgebras generated by $\{e_i\}_{i=1}^r$ and $\{ f_i \}_{i=1}^r$ will be denoted by $U_q(\n_+)$; set $U_q(\n_-)$, and $U_q(\b_\pm) := U_q(\n_\pm) U_q(\h)$. The algebras $U_q(\n_\pm)$ are graded by $\pm Q_+$ via 
\begin{equation}
    \label{Un-grad}
\deg e_i = \al_i = - \deg f_i.
\end{equation}
We denote by $T_i$, $i \in [1,r]$, Lusztig's braid group action on $U_q(\g)$ and on integrable $U_q(\g)$-modules, and set $T_u := T_{j_1} \cdots T_{j_m}$ for a reduced expression $u = s_{j_1} \cdots s_{j_m}$. For $\la \in P_+$, let $L(\la)$ be the irreducible highest weight $U_q(\g)$-module with highest weight $\la$ and highest weight vector $v_\la$. Unlike Sect. \ref{qg-forms}, $L(\la)$ is infinite dimensional unless $\g$ is finite dimensional, but for each $u \in W$ the weight space $L(\la)_{u\la}$ is still one dimensional. As in Sect. \ref{qg-forms}, denote $v_{u\la} := T_{u^{-1}}^{-1} v_\la$, and let $\xi_{u\la} \in (L(\la)_{u\la})^*$ be such that $\langle \xi_{u\la}, v_{u\la} \rangle = 1$.

Choose $d \in \Zset_+$ such that $(P^\vee, P^\vee) \subseteq \frac{1}{d} \Zset$. The Rosso--Tanisaki form is the Hopf pairing
\[
(\cdot, \cdot)_{RT} : U_q(\b_-) \times U_q(\b_+) \to \Cset(q^{1/d}),
\]
uniquely defined by
\[
(f_i, e_j)_{RT} := \frac{\delta_{ij}}{q_i^{-1} - q_i}, \quad 
(q^h, q^{h'})_{RT} := q^{-(h,h')}, \quad 
(f_i, q^h)_{RT} := (q^h, e_i)_{RT} = 0
\]
for $i, j \in [1,r]$ and $h, h' \in P^\vee$, 
see \cite[Eq. (5.5)]{GY2}.  

The divided power integral form $U_q(\n_+)_\cA$ is the $\cA$-subalgebra of $U_q(\n_+)$ generated by the divided powers $e_i^k / [k]_{q_i}!$ for $i \in [1,r]$ and $k \in \Zset_+$. Using the Rosso--Tanisaki form, we obtain the dual integral form of $U_q(\n_-)$:
\[
U_q(\n_-)^\vee_{\cA} := \{ x \in U_q(\n_-) \mid (x, y)_{RT} \in \cA \; \mbox{for all} \; y \in U_q(\n_+)_{\cA} \}.
\]
\subsection{Quantum Schubert cells and unipotent cells}
\label{unip-Schubert}
Fix $w \in W$ together with a reduced expression $w = s_{i_1} \cdots s_{i_N}$. Denote $w_{\leq k} := s_{i_1} \cdots s_{i_k}$ with the convention that $w_{\leq 0} := 1$. Denote the roots and root vectors 
\[
\be_k := w_{\leq k-1}(\al_{i_k})
\quad \mbox{and} \quad 
f_{\be_k} := T^{-1}_{w_{\leq k-1}^{-1}} (f_{i_k}) \in U_q(\n_-)_{\cA}
\quad \mbox{for} \quad k \in [1,N].
\]
The {\em{quantum Schubert cell algebra}} $U_q(\n_-(w))$, introduced by De Concini--Kac--Procesi \cite{DKP95} and Lusztig \cite{Lus93}, is the unital $\Cset(q)$-subalgebra of $U_q(\n_-)$ generated by $f_{\be_1}, \ldots, f_{\be_N}$; it does not depend on the choice of reduced expression for $w$. Intersecting it with the dual integral form of $U_q(\n_-)$ gives the $\cA$-algebra
\[
U_q(\n_-(w))^\vee_{\cA} := U_q(\n_-(w)) \cap U_q(\n_-)^\vee_{\cA}.
\]
Rescaling the root vectors $f_{\be_k}$ so that they pair to $1$ with $e_{\be_k}$ produces the {\em{dual PBW generators}}
\[
f'_{\be_k} := (q_{i_k}^{-1} - q_{i_k}) \, f_{\be_k} \in U_q(\n_-(w))^\vee_{\cA}, \quad k \in [1,N].
\]
Kimura \cite[Proposition 4.26 and Theorems 4.25, 4.27]{Kim17} proved that their ordered monomials form an $\cA$-basis of the dual integral form:
\begin{equation}
\label{unip-PBW}
U_q(\n_-(w))^\vee_{\cA} = \bigoplus_{m_1, \ldots, m_N \in \Nset} \cA \, (f'_{\be_1})^{m_1} \cdots (f'_{\be_N})^{m_N}.
\end{equation}

Following \cite{GLS}, consider the subalgebra $A_q(\n_+)$ of the full dual $U_q(\b_+)^*$ consisting of the functionals $\psi$ with the following two properties:
\begin{enumerate}
\item[(1)] $\psi(y q^h) = \psi(y)$ for all $y \in U_q(\n_+)$ and $h \in P^\vee$;
\item[(2)] there is a finite set $S \subseteq Q_+$ such that $\psi$ vanishes on $U_q(\n_+)_\gamma$ for every $\gamma \in Q_+ \backslash S$.
\end{enumerate}
Define the injective map $\iota : U_q(\n_-) \to U_q(\b_+)^*$ by 
\[
\langle \iota(x), y \rangle := (x, y)_{RT}, \quad x \in U_q(\n_-), \; y \in U_q(\b_+).
\]
Because the Rosso--Tanisaki form is a nondegenerate Hopf pairing, $\iota$ is an algebra isomorphism
\begin{equation}
\label{unip-iota-int}
\iota : U_q(\n_-) \stackrel{\cong}{\longrightarrow} A_q(\n_+).
\end{equation}

The {\em{quantum unipotent cell}} of Gei\ss--Leclerc--Schr\"oer associated to $w \in W$ is defined as the image of the quantum Schubert cell $U_q(\n_-(w))$ under $\iota$:
\[
\iota : U_q(\n_-(w)) \stackrel{\cong}{\longrightarrow} A_q(\n_+(w)) \subseteq A_q(\n_+).
\]
We will use the same notation for elements of $U_q(\n_-(w))$ and their images in $A_q(\n_+(w))$. 
The algebra $A_q(\n_+(w))$ is $Q_+$-graded by 
\begin{equation}
\label{unip-grading}
A_q(\n_+(w))_\gamma := \iota \big( U_q(\n_-(w))_{-\gamma} \big), \quad \gamma \in Q_+.
\end{equation}

The dual integral form of $U_q(\n_-(w))$ produces an $\cA$-form of the quantum unipotent cell,
\[
A_q(\n_+(w))_{\cA} := \iota \big( U_q(\n_-(w))^\vee_{\cA} \big),
\]
and $\iota$ restricts to an isomorphism of $\cA$-algebras
\[
\iota : U_q(\n_-(w))^\vee_{\cA} \stackrel{\cong}{\longrightarrow} A_q(\n_+(w))_{\cA}.
\]
The integral forms $U_q(\n_-)^\vee_{\cA}$, $U_q(\n_-(w))^\vee_{\cA}$ and $A_q(\n_+(w))_{\cA}$ are frequently defined through the bilinear forms of Kashiwara \cite{Kas91} or Lusztig \cite{Lus93} on $U_q(\n_-)$ rather than the Rosso--Tanisaki form; the resulting $\cA$-algebras are isomorphic \cite[Remark 5.3]{GY2}.

For $u \in W$ and $\mu \in P_+$, the {\em{unipotent quantum minor}} $D_{u\mu, w\mu} \in A_q(\n_+(w))_{(u-w)\mu}$ is defined, following \cite{GLS}, by 
\[
\langle D_{u\mu, w\mu}, y q^h \rangle := \langle \xi_{u\mu}, y v_{w\mu} \rangle, \quad \forall y \in U_q(\n_+), \; h \in P^\vee.
\]
It depends only on the weights $u\mu$ and $w\mu$, and not on the particular choice of $u$, $w$ and $\mu$ \cite[Sect. 9.3]{BerZe}. The minors $D_{\mu, w\mu}$ $q$-commute with all homogeneous elements for the $Q_+$-grading \eqref{unip-grading}, see \cite[Eq. (6.9)]{GY2}:
\begin{equation}
\label{unip-qcomm}
D_{\mu, w\mu} \, x = q^{((w+1)\mu, \gamma)} \, x \, D_{\mu, w\mu}, \quad 
\forall \mu \in P_+, \; x \in A_q(\n_+(w))_\gamma, \; \gamma \in Q_+.
\end{equation}

Let $\ep$ be a primitive $\ell$-th root of unity. The specialization of the {\em{quantum unipotent cell}} $A_q(\n_+(w))_{\cA}$ {\em{at the root of unity}} $\ep$ is
the $\Cset$-algebra
\[
A_\ep(\n_+(w)) := A_q(\n_+(w))_{\cA} / (q - \ep) A_q(\n_+(w))_{\cA}.
\]
As it is customary, by abuse of notation, we will denote by the same symbols the image of the elements of 
$A_q(\n_+(w))_{\cA}$ in $A_\ep(\n_+(w))$.
Combining Kimura's theorem \eqref{unip-PBW} with the isomorphism \eqref{unip-iota-int}, we see that the ordered monomials in the dual PBW generators $f'_{\beta_k}$ form a basis of $A_\ep(\n_+(w))$:
\[
A_\ep(\n_+(w)) = \bigoplus_{m_1, \ldots, m_N \in \Nset} \Cset \, (f'_{\be_1})^{m_1} \cdots (f'_{\be_N})^{m_N}.
\]

\bth{DKP-unip} {\em{(De Concini--Kac--Procesi \cite{DKP95})}} Let $\g$ be a symmetrizable Kac--Moody algebra, $w \in W$, and $\ep$ a primitive $\ell$-th root of unity, where $\ell$ is odd and coprime to the symmetrizing integers $d_i$ for all $i \in [1,r]$. Then
\[
C_\ep(\n_+(w)) := \bigoplus_{m_1, \ldots, m_N \in \Nset} \Cset \, (f'_{\be_1})^{\ell m_1} \cdots (f'_{\be_N})^{\ell m_N}
\]
is a central subalgebra of $A_\ep(\n_+(w))$
(called the Frobenius center of $A_\ep(\n_+(w))$).
\eth

De Concini, Kac and Procesi proved this theorem for finite dimensional $\g$ \cite{DKP95}; their argument carries over without change to arbitrary symmetrizable Kac--Moody algebras. In the general case, the theorem can also be derived from \cite[Proposition 4.4 and Theorem 8.5]{NTY2}, whose proof does not rely on \thref{DKP-unip}.

Recall that the support of $w \in W$ is defined to be the set 
\[
\SS(w) := \{ i \in [1,r] \mid s_i \leq w \}.
\]
Denote
\begin{equation}
\label{QSw}
Q_{\SS(w)}:= \bigoplus_{i \in \SS(w)} \Zset \al_i
\quad \mbox{and} \quad 
P_{\SS(w)}:= \bigoplus_{i \in \SS(w)} \Zset \vpi_i.
\end{equation}

\bth{Anw-clust} \cite[Theorems 8.1 and 8.5]{NTY2} For all 
symmetrizable Kac--Moody algebras $\g$, Weyl group elements $w \in W$ and primitive $\ell$-th roots of unity
$\ep$ such that $\ell$ is odd and coprime to $d_i$ for all $i \in [1,r]$, there exists a root of unity partially compactified quantum cluster algebra $\A_\ep(M_\ep, \B, \varnothing)$ (no frozen variables inverted) such that 
\[
A_\ep(\n_+(w)) = A_\ep(M_\ep, \B, \varnothing). 
\]
The quantum theoretic and cluster theoretic Frobenius centers coincide:
\[
C_\ep(\n_+(w)) = \C\A_\ep(M_\ep, \B, \varnothing).
\]
All cluster variables are homogeneous with respect to the $Q_+$-grading given by \eqref{Un-grad} and \eqref{unip-grading}. The frozen variables are the unipotent quantum minors 
\[
D_{\vpi_i, w \vpi_i} 
\quad \mbox{for} \quad 
i \in \SS(w).
\]
\eth
\subsection{Frobenius and ozone automorphism groups of quantum unipotent cells}
The abelian subgroup of $Q$ generated by 
\[
\{ \ga \in Q_+ \mid A_q(\n_+(w))_\gamma \neq 0 \} 
\]
equals 
\begin{equation}
\label{nonzer}
\sum_{j =1}^N \Zset \beta_j = \bigoplus_{i \in \SS(w)} \Zset \al_i = Q_{\SS(w)}.
\end{equation}
Here the last equality is the definition \eqref{QSw} and the first equality is proved in \cite[Lemma 3.2(ii)]{Y-spec} for finite dimensional Lie algebras $\g$; the proof given there works in general. 
This implies that we have a faithful action of $\Gset_m^{|\SS(w)|}$ on 
$A_\ep(\n_+(w))$ given by 
\begin{equation}
\label{t-act-Anw}
c \cdot u = \prod c_i^{m_i} u 
\quad \mbox{for} \quad 
c := (c_i)_{i \in \SS(w)} \in \Gset_m^{|\SS(w)|}(\Cset), \; 
u \in A_\ep(\n_+(w))_\gamma,
\end{equation}
where $\gamma = \sum_{i \in \SS(w)} m_i \alpha_i$. 

By the definition of the lattice $Q_{\SS(w)}$, 
\[
w^{-1} \big( Q_{\SS(w)} \big) = Q_{\SS(w)}. 
\]
Define the finite abelian group 
\begin{equation}
    \label{Qswl}
Q_{\SS(w),\ell}:=
(1+w^{-1})Q_{\SS(w)}/\big( (1+w^{-1})Q_{\SS(w)} \cap \ell Q_{\SS(w)} \big)
\end{equation}
and denote the projection
\[
(1+w^{-1})Q_{\SS(w)} \twoheadrightarrow Q_{\SS(w),\ell} \quad \mbox{by} \quad
\beta \mt \ol{\beta}.
\]
Recall that for a finite group $G$, $\wh{G}$ will denote its character group. We have an action of $\wh{Q_{\SS(w),\ell}}$ on $A_\ep(\n_+(w))$ by algebra automorphisms given by 
\[
\chi \cdot u = \chi \Big( \ol{(1+ w^{-1})(\gamma)} \Big) u, \quad \mbox{for} \quad 
\chi \in \wh{Q_{\SS(w),\ell}}, \; 
u \in A_\ep(\n_+(w))_\gamma, \; \gamma \in Q_{\SS(w)}. 
\]

\bth{AutAnw} Let $\g$ be a symmetrizable Kac--Moody algebra, $w \in W$, and $\ep$ a primitive $\ell$-th root of unity, where $\ell$ is odd and coprime to the symmetrizing integers $d_i$ for all $i \in [1,r]$.
\begin{enumerate}
    \item[(i)] The Frobenius automorphism group of the quantum unipotent cell $A_\ep(\n_+(w))$ is given by 
\[
\Aut(A_\ep(\n_+(w))/C_\ep(\n_+(w))) = C_\ell^{|\SS(w)|}
\]
with respect to the action \eqref{t-act-Anw}.

Every endomorphism of $A_\ep(\n_+(w))$ fixing $C_\ep(\n_+(w))$ is an automorphism.
    \item[(ii)] The ozone group of $A_\ep(\n_+(w))$ is
\[
\Oz(A_\ep(\n_+(w))) = \wh{Q_{\SS(w), \ell}},
\]
recall \eqref{Qswl}.
The elements of the group $\Oz(A_\ep(\n_+(w)))$ are conjugations by monomials 
in the normal elements $D_{\vpi_i, w \vpi_i}$ for $i \in \SS(w)$, 
recall \eqref{unip-qcomm}.
\end{enumerate}
\eth
\begin{proof} (i) \thref{Anw-clust} implies that 
the quantum group theoretic and cluster theoretic Frobenius centers of 
$A_\ep(\n_+(w)) = A_\ep(M_\ep, \B, \varnothing)$ coincide and that 
all cluster variables of $A_\ep(M_\ep, \B, \varnothing)$ are homogeneous with respect to the grading \eqref{unip-grading}. Since the action \eqref{t-act-Anw}
is faithful, we can apply \thref{torus} to this action to obtain
\[
\Aut(A_\ep(\n_+(w))/ C_\ep(\n_+(w))) = 
C_\ell^{|\SS(w)|}.
\]

(ii) \thref{normal-Oz} implies that every element of the ozone group 
$\Oz(A_\ep(\n_+(w)))$ is a conjugation by a monomial in the 
frozen variables $D_{\vpi_i, w \vpi_i}$, $i \in \SS(w)$,
of the cluster algebra $A_\ep(M_\ep, \B, \varnothing)$.
It follows from \eqref{unip-qcomm} that every such automorphism 
is of the form 
\[
u \mt \ep^{( (1+w) \la, \gamma )} u \quad \mbox{for} \quad 
u \in A_\ep(\n_+(w))_\gamma, \; \gamma \in Q_{\SS(w)}
\]
for some $\la \in P_{\SS(w)}$. Taking into account \eqref{nonzer} and 
$( (1+w) \la, \gamma ) = ( \la, (1 + w^{-1}) \gamma )$, 
we obtain that 
this group is isomorphic to the image of the homomorphism 
\[
\psi : P_{\SS(w)} \to \wh{Q_{\SS(w)}}, \quad 
\psi(\la)(\beta) := \ep^{( \la, \beta)}, \; \; 
\la \in P_{\SS(w)}, \beta \in (1 + w^{-1}) Q_{\SS(w)}. 
\]
The homomorphism 
\[
\pi : P_{\SS(w)} \to \wh{Q_{\SS(w), \ell}}, \quad 
\pi(\la)(\ol{\beta}) := \ep^{(\la, \beta )}, \; \; 
\la \in P_{\SS(w)}, \beta \in (1 + w^{-1}) Q_{\SS(w)}
\]
is well defined and surjective because $\ep$ is a primitive $\ell$-th root of unity, $\ell$ is coprime to $d_i$ for $i \in [1,r]$, 
$( \vpi_i, \al_j ) = d_j \langle \vpi_i, h_j \rangle = d_j \delta_{i,j}$ 
for $i, j \in [1,r]$ and the group $Q_{\SS(w), \ell}$ is finite with exponent that divides $\ell$. The canonical embedding $\wh{Q_{\SS(w),\ell}} \hra \wh{Q_{\SS(w)}}$, dual to the projection $Q_{\SS(w)} \twoheadrightarrow Q_{\SS(w),\ell}$, 
$\gamma \mt \ol{(1+w^{-1}) \gamma}$ gives rise to the commutative diagram
\[
\begin{tikzpicture}[scale=1.4]
\node (A) at (0,1) {$P_{\SS(w)}$}; 
\node (B) at (4,1) {$\wh{Q_{\SS(w)}}$};
\node (C) at (2,0) {$\wh{Q_{\SS(w),\ell}}$};
\draw
(A) edge[->] node[above]{$\psi$} (B)
(A) edge[->>] node[below]{$\pi$} (C)
(C) edge[right hook ->] (B);
\end{tikzpicture}
\]
Therefore the image of $\psi$ coincides with the image of the embedding 
$\wh{Q_{\SS(w),\ell}} \hra \wh{Q_{\SS(w)}}$, which proves the first statement in part (ii) of the theorem.
\end{proof}
We illustrate the second part of the theorem with the special case when $\g$ is a complex simple Lie algebra and $w=w_0$, the longest Weyl group element. In this case the algebra $A_\ep(\n_+)= A_\ep(\n_+(w_0))$ 
is isomorphic to the half quantum group $U_\ep(\n_+)$. 
\bth{complex} Let $\g$ be a complex simple Lie algebra, $\ep$ a primitive $\ell$-th root of unity for an odd positive integer $\ell$ that is not divisible by $3$ if $\g$ is of type $G_2$. Then 
\[
\Oz(A_\ep(\n_+)) \cong C_\ell^{r'},   
\]
where
\[
r' = 0
\]
if $\g$ is of type $B_r, C_r$, $D_r$ for $r$ even, $E_7$, $E_8$, $F_4$, $G_2$, 
\[
r' = \lfloor r/2 \rfloor
\]
if $\g$ is of type $A_r$, 
\[
r' = 1
\]
if $\g$ is of type $D_r$ for $r$ odd, and 
\[
r' = 2
\]
if $\g$ is of type $E_6$. 
\eth
\begin{proof} We have $A_\ep(\n_+) = A_\ep(\n_+(w_0))$. The map $i \mt i^*$, $i \in [1,r]$, where 
$\al_{i^*} = - w_0(\al_i)$ is an involution of the Dynkin graph which is the 
identity for types $B_r, C_r$, $D_r$ for $r$ even, $E_7$, $E_8$, $F_4$, $G_2$, the flip of the Dynkin graph in type $A_r$, the swap of the two extremal vertices while fixing the rest for type $D_r$ with $r$ odd and the map $1 \leftrightarrow 6$, 
$2 \leftrightarrow 5$ while fixing 3 and 4 in type $E_6$ in the Bourbaki labeling. 

Denote by $I_0$ the fixed point set of the involution $i \mt i^*$ and by $I_1$ a set of representatives of the 2-orbits of the involution. We claim that
\begin{equation}
\label{Qswii*}
Q_{\SS(w_0), \ell} = \bigoplus_{i \in I_1} (\Zset/\ell) (\al_i - \al_{i^*}).
\end{equation}
The statement of the theorem follows from this equality by applying 
\thref{AutAnw}(ii) and taking into account that, for all types, 
$r'=|I_1|$.

We have 
\[
(1 + w_0^{-1}) Q = \bigoplus_{i \in I_1} \Zset (\al_i - \al_{i^*})
\]
and 
\[
(1 + w_0^{-1}) Q \cap \ell Q = \bigoplus_{i \in I_1} \ell \Zset (\al_i - \al_{i^*}).
\]
The two equalities imply \eqref{Qswii*}, which completes the proof of the theorem. 
\end{proof}


\end{document}